\documentclass{article}
\usepackage{hyperref}
\usepackage{amssymb}
\usepackage{amsthm}
\usepackage{mathtools}
\usepackage{tikz-cd}
\usepackage{bbm}
\usepackage[numbers,sort]{natbib}
\usepackage{todonotes}

\newtheorem{theorem}{Theorem}
\newtheorem{lemma}[theorem]{Lemma}
\newtheorem{corollary}[theorem]{Corollary}
\newtheorem{proposition}[theorem]{Proposition}
\theoremstyle{definition}
\newtheorem{remark}[theorem]{Remark}
\newtheorem{definition}[theorem]{Definition}

\newcommand{\clA}{\mathcal{A}}
\newcommand{\clC}{\mathcal{C}}

\newcommand{\clI}{\mathcal{I}}
\newcommand{\clD}{\mathcal{D}}
\newcommand{\clT}{\mathcal{T}}
\newcommand{\clL}{\mathcal{L}}
\newcommand{\clR}{\mathcal{R}}
\newcommand{\bbD}{\mathbb{D}}

\newcommand{\bbN}{\mathbb{N}}
\newcommand{\bbX}{\mathbb{X}}
\newcommand{\bbY}{\mathbb{Y}}
\newcommand{\bbG}{\mathbb{G}}
\newcommand{\bbM}{\mathbb{M}}
\newcommand{\one}{\mathbbm{1}}

\newcommand{\op}[1]{#1^{\textrm{op}}}
\DeclareMathOperator{\Bat}{Bat}
\DeclareMathOperator{\Grid}{Grid}
\DeclareMathOperator{\Pos}{Pos}
\DeclareMathOperator{\Set}{Set}
\DeclareMathOperator{\Ord}{Ord}
\DeclareMathOperator{\Comp}{Comp}

\DeclareMathOperator{\Sig}{Sig}
\DeclareMathOperator{\Alg}{Alg}
\DeclareMathOperator{\Plex}{Plex}
\DeclareMathOperator{\Pplex}{Pplex}
\DeclareMathOperator{\Term}{Term}
\DeclareMathOperator{\Type}{Type}
\DeclareMathOperator{\Cptd}{Cptd}
\DeclareMathOperator{\Free}{Free}
\DeclareMathOperator{\Und}{Und}
\DeclareMathOperator{\F}{F}
\DeclareMathOperator{\U}{U}
\DeclareMathOperator{\M}{M}
\DeclareMathOperator{\K}{K}
\DeclareMathOperator{\Cof}{Cof}
\DeclareMathOperator*{\colim}{colim}
\DeclareMathOperator{\sk}{sk}
\DeclareMathOperator{\tr}{tr}
\DeclareMathOperator{\ob}{ob}
\DeclareMathOperator{\id}{id}
\DeclareMathOperator{\inc}{inc}
\DeclareMathOperator{\var}{var}
\DeclareMathOperator{\dpth}{dpth}
\DeclareMathOperator{\sort}{sort}
\DeclareMathOperator{\supp}{supp}
\DeclareMathOperator{\face}{face}
\DeclareMathOperator{\filler}{fill}
\DeclareMathOperator{\Kan}{Kan}
\DeclareMathOperator{\cat}{cat}
\DeclareMathOperator{\coh}{coh}
\DeclareMathOperator{\grpd}{grpd}
\DeclareMathOperator{\mcat}{mcat}
\DeclareMathOperator{\ar}{ar}

\newcommand{\Compv}{\Comp^{\textrm{var}}}
\newcommand{\Termv}{\Term^{\textrm{var}}}
\newcommand{\Cptdv}{\Cptd^{\textrm{var}}}
\DeclarePairedDelimiter{\abs}{\lvert}{\rvert}

\begin{document}

\title{Computads for generalised signatures}
\author{Ioannis Markakis\footnote{University of Cambridge,
  \texttt{im496@cam.ac.uk}}}
\date{March 21, 2023}


\maketitle

\begin{abstract}
    We introduce a notion of signature whose sorts form a direct
    category, and study computads for such signatures. Algebras for
    such a signature are presheaves with an interpretation of every function
    symbol of the signature, and we describe how computads give rise to
    signatures. Generalising work of Batanin, we show that computads with
    certain generator-preserving morphisms form a presheaf category, and
    describe a forgetful functor from algebras to computads. Algebras free on a
    computad turn out to be the cofibrant objects for certain cofibrantly
    generated factorisation system, and the adjunction above induces the
    universal cofibrant replacement, in the sense of Garner, for this
    factorisation system. Finally, we conclude by explaining how
    many-sorted structures, weak $\omega$\nobreakdash-categories, and algebraic
    semi-simplicial Kan complexes are algebras of such signatures,
    and we propose a notion of weak multiple category.
\end{abstract}

\section{Introduction}\label{sec:introduction}

An important question for any kind of mathematical structure is determining
the underlying data that it can be freely built from. For most algebraic
structures, such as groups or rings, that data is given by a set of generators,
while for \emph{higher dimensional} structures, more structured data is needed,
for example categories can be freely built from a directed graph. In all those
examples, the structures of interest are algebras of some monad $\M$ on some
category $\clC$ and the free generation is expressed via the Eilenberg-Moore
adjunction:
\[\begin{tikzcd}
	\clC & \Alg_{\M}
	\arrow[""{name=0, anchor=center, inner sep=0}, "{\F^{\M}}",
  shift left=2, from=1-1, to=1-2]
	\arrow[""{name=1, anchor=center, inner sep=0}, "{\U^{\M}}",
  shift left=2, from=1-2, to=1-1]
	\arrow["\dashv"{anchor=center, rotate=-90}, draw=none, from=0, to=1]
\end{tikzcd}\]
A similar adjunction exists for globular higher categories, where the category
$\clC$ is the category of globular sets and $\M$ some finitary monad. However,
in that case, there is also a more general notion of generating datum, called
\emph{computads} or \emph{polygraphs}~\cite{street_limits_1976,
burroni_higherdimensional_1993,batanin_computads_1998},
build recursively by gluing disks along their boundary spheres, similar to CW
complexes.

In this paper, we introduce and study computads for monads over more general
categories $\clC$. More precisely, we take $\clC = [\op\clI, \Set]$ to be a
category of presheaves over an arbitrary small direct category $\clI$, the
objects of which we call \emph{sorts} and the morphisms of which we call
\emph{face maps}. Each sort $i$ gives rise to an inclusion
$\partial \bbD^i\subseteq\bbD^i$ in $\clC$ where $\bbD^i = \clI(-,i)$ is the
representable presheaf by $i$ and $\partial \bbD^i$ the sub-presheaf obtained
by removing the top-dimensional element of $\bbD^i$. This broader framework
allows us to treat computads for higher categories uniformly, regardless of the
underlying geometry being globular, semi-simplicial or cubical for example.

Furthermore, we study what data can generate a monad $\M$ over such a
category $\clC$. When $\clI$ is the terminal category, so $\clC$ is the
category of sets, an answer to this questions is the notion of
\emph{signature} from universal algebra~\cite{hodges_model_1993}. Roughly that
amounts to a set $\Sigma$ of function symbols together with an \emph{arity}
function $\ar : \Sigma\to \bbN$ assigning to each symbol its number of
inputs. A \emph{$\Sigma$\nobreakdash-algebra} $\bbX$ amounts to a set $X$ together with a
function
\[ f^\bbX : X^{\ar(f)} \to X\]
for every function symbol $f\in \Sigma$. The monad $\M_\Sigma$ generated by
$\Sigma$ is the one induced by the free $\Sigma$\nobreakdash-algebra adjunction
$\F_\Sigma \dashv \U_\Sigma$.

A vast generalisation of signature has been proposed
where $\clC$ is a locally presentable enriched category, and the natural numbers
are replaced by a small dense subcategory $\clA$ of arities~\cite{bourke_monads_2019}. Unwrapping the
definition in the case that $\clC = [\op\clI, \Set]$ is a category of
presheaves, we see that a signature $\Sigma$ amounts to
a presheaf of function symbols for every $A\in \clA$. Equivalently, it is a
presheaf $\Sigma : \op\clI\to \Set$ together with an arity function
$\ar : \Sigma_i\to \ob\clA$ for every sort $i$ such that
$\ar(\delta^*f) = \ar(f)$
for every function symbol $f\in \Sigma_i$ and face map $\delta : j\to i$. A
$\Sigma$\nobreakdash-algebra $\bbX$ is a presheaf $X$ together with a function
\[ f^\bbX : \clC(\ar(f), X) \to X_i \]
compatible with the face maps in that
\[ \delta^* f^\bbX = (\delta^*f)^\bbX. \]
The monad $\M_\Sigma$ generated by $\Sigma$ is again the one induced by the
free $\Sigma$\nobreakdash-algebra adjunction $\F_\Sigma \dashv \U_\Sigma$. In this setting,
every element $t\in (\M_\Sigma A)_i$ can be seen as a formal composite of
function symbols with arity $A\in \clA$, since it gives rise to a function
\[
  t^\bbX : \clC(A,X) \xrightarrow{\sim}
  \Alg_\Sigma(\F_\Sigma A, X) \to
  X_i
\]
by evaluation.

Here we propose a generalised notion of signature in the case that $\clI$ is a
small direct category. Our notion of signature generalises the previous one, in
the same way that computads generalise presheaves on $\clI$ as generating data.
Our signatures consist again of a set of function symbols $\Sigma_i$ for every
sort $i$ together with an arity function $\ar : \Sigma_i\to \ob\clC$. However,
instead of forming a presheaf, the function induced by a face map
$\delta : j\to i$ sends a function symbol $f\in \Sigma_i$ of arity $A$ to a
formal composite of function symbols $t_{f,\delta}\in \M_{\Sigma,j}A$ of lower
dimension. In order to make this definition precise, we work recursively on the
dimension function $\dim : \ob\clI\to \Ord$ of the direct category, and define
signatures, $\Sigma$\nobreakdash-computads and an adjunction inducing the monad $\M_\Sigma$
simultaneously by mutual induction.

A benefit of working with monads generated by a signature is that computads,
their morphisms and the free algebras on them can be described as
\emph{inductively generated sets}, that is they are themselves freely generated
by a set of constructors. To illustrate this concept, let $\Sigma$ be a signature
in the sense of universal algebra. Computads for such a signature $\Sigma$ are
just sets and the monad $\M_\Sigma$ is the one sending a set $X$ to its set of
\emph{$\Sigma$\nobreakdash-terms}. This is a set inductively generated by the constructors
\begin{itemize}
  \item  there exists a term $\var(x)$ for every $x\in C$,
  \item there exists a term $f[t_1,\dots,t_{\ar(f)}]$ for every
    $f\in \Sigma$ and terms $t_1,\dots, t_{\ar(f)}$.
\end{itemize}
What that means is that the set $\M_\Sigma X$ is an initial algebra for the
polynomial endofunctor
\[
  G(X) = C\amalg \coprod_{f\in \Sigma} X^{\ar(f)}
\]
specified by the constructors above. To construct the initial algebra, we let
first $X_{-1} = \emptyset$ the initial set and form an increasing family of sets
by letting $X_{n} = G(X_{n-1})$ for every natural number $n$. The union of this
family is the initial algebra~\cite{adamek_free_1974}. A more detailed
introduction to inductively generated sets can be found in our previous work~\cite[Appendix~A]{dean_computads_2022}.

This presentation of computads and free algebras allows us to give simple proofs
of existing results on globular computads~\cite{batanin_computads_1998, garner_homomorphisms_2010}. In Section~\ref{sec:computads-are-presheaves}, we show
that the category of computads and certain generator-preserving maps is a
presheaf topos for every generalized signature $\Sigma$. We do that by showing
first that the functor of \emph{terms}, sending a computad to the underlying
presheaf of the free algebra on it, is familially representable. Following~\cite{henry_nonunital_2019}, we call the computads representing this functor
\emph{polyplexes} and we identify a subset of them, the \emph{plexes}, that
correspond to representable computads. Our proofs on this section are
similar to the ones in our previous work~\cite{dean_computads_2022}. This result
crucially uses that the monads we consider are free on a signature, since it
fails for example for the free strict $\omega$\nobreakdash-category monad~\cite{makkai_category_2008}.

In Section~\ref{sec:computads-are-cofibrant} we study the weak
factorisation system on the category of algebra generated by the boundary
inclusions $\partial \bbD^i\subseteq \bbD^i$. We call morphisms in its left
class \emph{cofibrations} and morphisms in the right class \emph{trivial
fibrations}, due to their role in the examples in Section~\ref{sec:examples}.
We show that free algebras on a computad and cofibrant
algebras coincide, as in the case of strict
$\omega$\nobreakdash-categories~\cite{metayer_cofibrant_2008}. We also construct an adjunction between the
category of algebras and the topos $\Compv_\Sigma$ of computads and
generator-preserving maps, extending the free algebra adjunction:
\[
  \begin{tikzcd}
    & {\Alg_\Sigma} \\
    \clC && \Compv_\Sigma
    \arrow[""{name=0, anchor=center, inner sep=0}, "{\U_\Sigma}",
    shift left=1, from=1-2, to=2-1]
    \arrow[""{name=1, anchor=center, inner sep=0}, "{\F_\Sigma}",
    shift left=2, from=2-1, to=1-2]
    \arrow[""{name=2, anchor=center, inner sep=0}, "{\Und_\Sigma}"',
    shift right=1, from=1-2, to=2-3]
    \arrow[""{name=3, anchor=center, inner sep=0}, "{\Free_\Sigma}"',
    shift right=2, from=2-3, to=1-2]
    \arrow[hook, from=2-1, to=2-3]
    \arrow[no head, maps to, shorten <= 1pt,shorten >= 1pt, from=2, to=3]
    \arrow[no head, maps to, shorten <= 1pt,shorten >= 1pt, from=0, to=1]
  \end{tikzcd}
\]
We show that the induced comonad is the
\emph{universal cofibrant replacement} for this weak factorisation
system~\cite{garner_homomorphisms_2010}.

Finally, in Section~\ref{sec:examples}, we descirbe certain
signatures of interest and their algebras. We see first that many-sorted
algebraic signatures are $S$-sorted signatures for $S$ a discrete category. We
then show that $\omega$\nobreakdash-categories and algebraic semi-simplicial Kan
complexes are algebras for signatures with sorts from the categories of globes
and semi-simplices respectively. We conclude by proposing a new notion of fully
weak unbiased multiple category as algebras of a certain signature with sorts in
a cube category.

\subsection*{Related Work}

Our \emph{$\clI$-sorted signatures} are special cases of Fiore's
\emph{$\Sigma_0$\nobreakdash-models with substtitution}~\cite{fiore_secondorder_2008}, as
explained by Subramaniam~\cite{subramaniam_dependent_2021}, who compares Fiore's
models with his \emph{$\clI$\nobreakdash-sorted theories}. In particular, $\clI$\nobreakdash-sorted
signatures seem to be closely related to the
\emph{dependently typed term signatures} of Subramaniam.

In order to give a syntactic presentation of theories over certain contextual
categories, Subramaniam introduced in his
thesis~\cite[Section~1.6]{subramaniam_dependent_2021}, \emph{dependently typed
type signatures} which are extensions of structural Martin-L\"of dependent
type theory by a sequence of type declarations $\Gamma \vdash A$. Such
signatures can be seen as presentations of locally finite direct categories,
where the types $A$ being declared correspond to objects of the category, and
the context $\Gamma$ over which they are declared, describe morphisms into this
object; more precisely the subpresheaf of the representable presheaf obtained
by removing the identity. A \emph{dependently typed term signature} over a type
signature is further extension by a sequence of \emph{term declarations} which
are axioms of the form \[\Delta \vdash t : A[\sigma],\]
where $\Delta$ is a context over the type signature, $(\Gamma\vdash A)$ is a
type declaration and $\sigma :\Delta \to \Gamma$ a morphism of contexts.

An $\clI$-sorted signature $\Sigma$ in our sense, where the category $\clI$ is
locally finite and the arities of its function symbols are finite presheaves
gives rise to a dependently typed term signature, whose term declarations are
the function symbols of $\Sigma$. The context $\Delta$ of a term declaration is
given by the arity of the function symbol, the type $A$ is give by the output
sort and the substitution $\sigma$ is assembled from the faces of the function
symbol. Conversely, every dependently typed term signature, in which the
dimension of the context $\Delta$ of every term declaration is bounded by the
dimension of the type $A$, seen as a presheaf and an object of the category
$\clI$ respectively, gives rise to an $\clI$-sorted signature. This natural
condition on the dimension is necessary for monads free on generalised
signatures to commute with the truncation functors in a sense made precise
below.

\subsection*{Notation}

We will denote by $\clI$ a small \emph{direct category}, whose objects we will
call \emph{sorts} and whose morphisms we will call \emph{face maps}. By direct
category, we mean that it is equipped with a function $\dim : \ob\clI\to \Ord$
to the class of ordinals satisfying that $\dim j < \dim i$ when there exists a
non-identity face map $j\to i$. It is easy to see that $\clI$ must be skeletal
and have no non-identity endomorphisms. Given a presheaf
$X : \op{\clI}\to \Set$, we will denote by $X_i$ the value of $X$ at a sort $i$
and by $\delta^* : X_i\to X_j$ its value at a face map $\delta : j\to i$.

Given an ordinal $\alpha$, we will let $\clI_\alpha$ denote the full subcategory
of $\clI$ whose objects have dimension at most $\alpha$. Pulling back along the
subcategory inclusion induces a \emph{truncation} functor
\[
  \tr_\alpha : [\op{\clI},\Set] \to [\op{\clI_\alpha},\Set].
\]
Left and right Kan extensions along the inclusion, define a left and right
adjoint respectively to the truncation functor called the \emph{skeleton}
and \emph{coskeleton} functors
respectively~\cite[Chapter~1]{riehl_categorical_2014}. The subcategory inclusion
$\clI_\beta\subseteq \clI_\alpha$ for $\beta \le \alpha$ defines similarly an
adjoint triple and it is easy to see that the following conditions are satisfied
\begin{align*}
  \tr_\beta &= \tr^\alpha_\beta \tr_\alpha &
  \tr_\beta\sk_\beta &= \id &
  \sk_\beta &= \sk_\alpha \sk^\alpha_\beta,
\end{align*}
which we will call the \emph{cocycle conditions}.

\subsection*{Acknowledgements}
I would like to thank my supervisor Prof. Jamie Vicary for his support during
this project. Furthermore, I would like to acknowledge funding from the Onassis
foundation - Scholarship ID: F~ZQ~039-1/2020-2021. I would finally like to
thank Chaitanya Leena Subramaniam for helpful conversations.

\section{Signatures and computads}\label{sec:definitions}

In this section, we will define the class of \emph{$\clI$-sorted signatures}.
We will simultaneously define for every signature $\Sigma$, a category of
\emph{computads}, extending the category of
presheaves on $\clI$, and the presheaf of terms of a
computad. The definition is given by transfinite recursion. More
precisely, for every ordinal $\alpha$, we will define a class of
\emph{$\clI$-sorted signatures of dimension $\alpha$} together with
restriction functions
\[(-)_\beta : \Sig_\clI(\alpha) \to \Sig_\clI(\beta)\]
for $\beta \le \alpha$. Moreover, for every such signature $\Sigma$, we
will define a category of computads together with truncation functors
\[ \tr^\Sigma_\beta : \Comp_\Sigma \to \Comp_{\Sigma_\beta} \]
for $\beta\le \alpha$. We will also define with an adjunction
\begin{align*}
  \Cptd_\Sigma &: [\op{\clI_\alpha},\Set] \to \Comp_\Sigma &
  \eta_\Sigma  &: \id\Rightarrow \Term_\Sigma\Cptd_\Sigma \\
  \Term_\Sigma &: \Comp_\Sigma\to [\op{\clI_\alpha},\Set] &
  \varepsilon_\Sigma &: \Cptd_\Sigma\Term_\Sigma \Rightarrow \id
\end{align*}
commuting with the truncation functors, in the sense that
\begin{align*}
  \tr_\beta^\Sigma \Cptd_\Sigma
    &= \Cptd_{\Sigma_\beta}\tr_\beta^\alpha &
  \tr_\beta^\alpha \eta_\Sigma
    &= \eta_{\Sigma_\beta}\tr_\beta^\Sigma \\
  \tr_\beta^\alpha \Term_\Sigma
    &= \Term_{\Sigma_\beta}\tr_\beta^\Sigma &
  \tr_\beta^\Sigma\varepsilon_\Sigma
    &=\varepsilon_{\Sigma_\beta}\tr_\beta^\Sigma.
\end{align*}
We will denote the monad induced by this adjunction by
$(\M_\Sigma,\eta_\Sigma,\mu_\Sigma)$. The restriction
functions and the truncation functors will be shown to satisfy the following
\emph{cocycle conditions}
\begin{align*}
  (-)^\alpha_\gamma &= (-)^\beta_\gamma(-)^\alpha_\beta &
  (-)^\alpha_\alpha &= \id \\
  \tr^\Sigma_\gamma &= \tr^{\Sigma_\beta}_\gamma \tr^\Sigma_\beta &
  \tr^\Sigma_\alpha &= \id
\end{align*}
for $\gamma\le\beta\le\alpha$.

For the rest of the section, we fix an ordinal $\alpha$, and suppose that the
data above has been defined for every $\beta\le\alpha$ and it satisfies the
cocycle conditions. We proceed to define that data for $\alpha$ as well. The
apparent circularity of our definition will be explained in Remark~\ref{rmk:recursive-depth}.

\paragraph{Signatures}
An \emph{$\clI$-sorted signature $\Sigma$} of dimension $\alpha$ consists of
\begin{itemize}
  \item a signature $\Sigma_\beta$ of dimension $\beta$
    for every $\beta < \alpha$ satisfying for $\gamma \le \beta$ that
    \[ (\Sigma_\beta)_\gamma = \Sigma_\gamma, \]
  \item a set $\Sigma_i$ of function symbols for every sort $i$ of
    dimension $\alpha$, and for every function symbol $f\in \Sigma_i$,
    \begin{itemize}
      \item a presheaf $B_f$ on $\clI_\alpha$, called the \emph{arity} of the
        function symbol,
      \item a \emph{boundary} term $t_{f,\delta}\in
        \M_{\Sigma_{\beta},j}(\tr_\beta B_f)$ for every non-identity face map
        $\delta : j\to i$, where $\beta = \dim j$ satisfying for every face
        map $\delta' : k\to j$ that
        \[ (\delta')^*(t_{f,\delta'}) = t_{f,\delta\delta'}. \]
    \end{itemize}
\end{itemize}
The restriction function $(-)_\beta$ for $\beta \le \alpha$ is the obvious
projection.

\paragraph{Computads}
A \emph{$\Sigma$\nobreakdash-computad $C$} consists of
\begin{itemize}
  \item a $\Sigma_\beta$-computad $C_\beta$ for every $\beta <\alpha$ satisfying
    for $\gamma \le \beta$ that
    \[\tr_\gamma^{\Sigma_\beta}(C_\beta) = C_\gamma,\]
  \item a set $V_i^C$ of \emph{generators} for every sort $i$ of
    dimension $\alpha$,
  \item a \emph{gluing} function
    $\phi_\delta^C : V_i^C\to \Term_{\Sigma_{\beta},j}(C_{\beta})$ for every
    non-identity face map $\delta : j\to i$, where $\beta = \dim j$, satisfying
    for every face map $\delta' : k\to j$ that
    \[ (\delta')^*\phi_{\delta}^C = \phi_{\delta\delta'}^C. \]
\end{itemize}
A \emph{morphism of computads} $\sigma : C\to D$ consists similarly of
\begin{itemize}
  \item morphisms $\sigma_\beta : C_\beta \to D_\beta$ for every
    $\beta < \alpha$ satisfying for $\gamma \le \beta$ that
    \[ \tr_\gamma^{\Sigma_\beta}(\sigma_\beta) = \sigma_\gamma \]
  \item functions $\sigma_i : V_i^C\to \Term_{\Sigma,i}(D)$ for every sort $i$
    of dimension $\alpha$ satisfying for every non-identity face map
    $\delta : j\to i$ that
    \[ \delta^* \sigma_i =
       \Term_{\Sigma_\beta,j}(\sigma_\beta) \circ \phi_\delta^C
    \]
    where $\beta = \dim j$.
\end{itemize}
The composition of a pair of morphisms $\sigma : C\to D$ and $\tau : D\to E$
is given by
\begin{align*}
  (\tau\circ\sigma)_\beta &= \tau_\beta\circ\sigma_\beta &
  (\tau\circ\sigma)_i &= \Term_{\Sigma,i}(\tau)\circ\sigma_i
\end{align*}
for $\beta < \alpha$ and $i$ of sort $\alpha$. The identity of a computad $C$ is
given by the identities $\id_{C_\beta}$ and the inclusion
$\var : V_i^C \to \Term_{\Sigma,i}(C)$ of generators into terms defined below.
The truncation functor $\tr^\Sigma_\beta$ for $\beta \le \alpha$ is the obvious
projection.

\paragraph{The inclusion functor}
The computad $\Cptd_\Sigma X$ associated to a presheaf
$X$ comprises of the computads $\Cptd_{\Sigma_\beta}\tr_\beta X$ for every
$\beta < \alpha$, the sets $X_i$ for every sort $i$ of dimension $\alpha$, and
the gluing functions
\[
  \phi_\delta^{\Cptd_\Sigma X} : X_i \xrightarrow{\delta^*} X_j
  \xrightarrow{\eta_{\Sigma_{\beta},\tr_\beta X,j}}
  \M_{\Sigma_\beta,i}(\tr_\beta X)
\]
for every face map $\delta : j\to i$, where $\beta = \dim j$. The morphism of
computads induced by a morphism of presheaves $\sigma : X\to Y$
comprises similarly of $\Cptd_{\Sigma_\beta}\tr_\beta\sigma$ for every
$\beta < \alpha$, and the composites
\[
  (\Cptd_\Sigma \sigma)_i :
  X_i \xrightarrow{\sigma_i} Y_i
  \xrightarrow{\var}
  \M_{\Sigma,i}Y
\]
where $\var$ the inclusion of generators into terms defined below.

\paragraph{The presheaf of terms}
The presheaf $\Term_\Sigma(C)$ of terms of a computad $C$ is defined
as follows. For every sort $j$ of dimension $\beta < \alpha$, we let
\[
  \Term_{\Sigma,j}(C) = \Term_{\Sigma_\beta,j}(C_\beta).
\]
The \emph{boundary function} $\delta^*$ induced by a face map $\delta : k\to j$
is given by
\[
  \Term_{\Sigma_\beta,j}(C_\beta)
  \xrightarrow{\delta^*}
  \Term_{\Sigma_\beta,k}(C_\beta) =
  \Term_{\Sigma_\gamma,k}(C_\gamma)
\]
where $\gamma = \dim k$. For a sort $i$ of dimension $\alpha$, the set
$\Term_{\Sigma,i}(C)$ is defined inductively by the constructors
\begin{itemize}
  \item there exists a term $\var v$ for every generator $v\in V_i^C$,
  \item there exists a term $f[\tau]$ for every function symbol
    $f\in \Sigma_i$ and morphism of computads $\tau : \Cptd_\Sigma B_f\to C$.
\end{itemize}
The function $\delta^*$ induced by a non-identity face map
$\delta : j\to i$ is defined recursively by
\begin{align*}
  \delta^*(\var v)
    &= \phi^C_\delta(v) &
  \delta^*(f[\tau])
    &= \Term_{\Sigma_{\beta},j}(\tau_{\beta})(t_{f,\delta})
\end{align*}
where $\beta = \dim j$.
The cocycle conditions ensure that this assignment defines a presheaf on
$\clI_\alpha$. The action of a morphism of computads $\sigma : C\to D$ on a
term $t$ of sort $j$ of dimension $\beta < \alpha$ is given by
\[\Term_{\Sigma}(\sigma)(t) = \Term_{\Sigma_\beta}(\sigma_\beta)(t),\]
while its action on a term of sort $i$ of dimension $\alpha$ is defined
recursively by
\begin{align*}
  \Term_{\Sigma}(\sigma)(\var v) &= \sigma_i(v) &
  \Term_{\Sigma}(\sigma)(f[\tau]) &= f[\sigma\circ\tau].
\end{align*}

\paragraph{The adjunction}
It remains to define the unit and counit of the term adjunction. For
the former, given a presheaf $X$ on $\clI_\alpha$, let
\[\eta_{\Sigma,X} : X\to \Term_{\Sigma}\Cptd_{\Sigma}X\]
the morphism of presheaves defined on some $x\in X_i$ by
\[
  \eta_{\Sigma,X,i}(x) =
    \begin{cases}
      \eta_{\Sigma_{\beta}, \tr_{\beta}X,i}(x),
        &\text{for }\dim i = \beta < \alpha \\
      \var(x), &\text{for } \dim i = \alpha.
    \end{cases}
\]
For the latter, given a computad $C$ for $\Sigma$, let
\[\varepsilon_{\Sigma,C} : \Cptd_\Sigma\Term_\Sigma C \to C\]
the morphism of computads consisting of $\varepsilon_{\Sigma_\beta,C_\beta}$
for every $\beta < \alpha$ and the identity of the set
\[
  V^{\Cptd_\Sigma\Term_\Sigma C}_i = \Term_{\Sigma,i}(C)
\]
for $i$ of dimension $\alpha$. Naturality of the unit and counit of the
adjunction as well as the snake equations can be easily checked. This concludes
the recursive definition.

\begin{remark}\label{rmk:recursive-depth}
Albeit the apparent circularity of the definition of the category of
computads and the functor of terms for a fixed signature $\Sigma$,
this is a valid mathematical definition as explained below. Assuming that all
the data has been defined for every $\beta < \alpha$,
the notion of signature of dimension $\alpha$ can be defined. Fixing such a
signature $\Sigma$, computads and their terms of sort $\beta < \alpha$ may be
defined as well as the inclusion functor on objects. Then fixing a computad $C$,
its terms of dimension $\alpha$, the boundary functions, and
morphisms from an arity to $C$ are defined using induction recursion~\cite{hancock_small_2013}.

To understand the inductive recursive definition, we attach an ordinal to
every term of $C$ of dimension $\alpha$ and every morphism of computads
with target $C$, its \emph{depth}, recursively by
\begin{align*}
  \dpth(\var v)
    &= 0 \\
  \dpth(f[\tau])
    &= \dpth\tau + 1 \\
  \dpth(\sigma : D\to C)
    &= \sup\{ \dpth(\sigma_i(v)) : \dim i = \alpha, v\in V_i^C \}.
\end{align*}
The definition of terms and morphisms amounts then to constructing by
transfinite recursion on an ordinal $\gamma$ two increasing families of sets,
and forming their unions. The first family consists of the sets
$\Term_{\Sigma,i}^\gamma(C;(t_\delta))$ of terms of sort $i$, depth at most
$\gamma$ and boundaries given by the terms $t_\delta$. The second consists of
the sets $\Comp_\Sigma^\gamma(B_f,C;(\sigma_\beta))$ of morphisms of depth at
most $\gamma$ from an arity $\Cptd_\Sigma B_f$ to $C$ and truncations given by
the morphisms $\sigma_\beta$. Those sets are defined by the recursive formulae
\begin{align*}
  \Term_{\Sigma,i}^\gamma(C;(t_\delta)) &=
    \left(\bigcap_{\delta}(\phi_\delta^C)^{-1}(t_\delta)\right) \amalg
    \coprod_{\substack{f\in \Sigma_i \\ (\sigma_\beta)}}
    \bigcup_{\gamma'<\gamma} \Comp_{\Sigma}^{\gamma'}(B_f,C;(\sigma_\beta)) \\
  \Comp_\Sigma^\gamma(B_f,C;(\sigma_\beta)) &=
    \prod_{\substack{\dim j = \alpha \\ b\in B_{f,j}}}
    \Term_{\Sigma,j}^\gamma(C;(\Term(\sigma_{\beta})(t_{f,\delta})))
\end{align*}
where the first coproduct is over families of morphisms
$\sigma_\beta : \Cptd_{\Sigma_\beta}\tr_\beta B_f\to C_\beta$ for
$\beta < \alpha$ satisfying the cocycle conditions and that
$\Term_\Sigma(\sigma_{\beta})(t_{f,\delta}) = t_\delta$ for every non-identity
face map $\delta : j\to i$, where $\beta = \dim j$. The union of the first
family over all $\gamma$ and all families of boundary terms $(t_\delta)$ is the
set of terms of sort $i$ of $C$. The boundaries of those terms are
given by the obvious projections. Similarly, the union of the second family
gives morphisms from $\Cptd_\Sigma B_f$ to $C$ and the obvious projections give
the truncations of those morphisms. Once terms are defined, we may
also define morphisms between arbitrary computads, and the truncation functor in
general.

Observe that we defined the sets of terms and the sets of morphisms as
union of increasing families over all ordinals. Such unions produce proper
classes instead of sets unless they are eventually stationary. To see that this
is the case in our definition, let $\lambda$ a regular cardinal strictly greater
than the cardinality of the set $\coprod_{\dim j = \alpha} B_{f,j}$ for every
function symbol of dimension $\alpha$. Existence of such cardinal follows by
$\clI$ being small and the $\Sigma_i$ being sets. A simple inductive
argument shows that all terms and morphisms from an arity have depth
strictly less than $\lambda$, and arbitrary morphisms have depth at most
$\lambda$.

Once those are defined, we may define composition of morphisms and the action of
morphisms on terms mutually recursively, or equivalently by induction on depth.
By induction on depth, we can then prove that composition is associative and
unital, and that the action on terms is functorial. Once those properties are
established, it is easy to define the action of the inclusion functor on
morphisms, as well as the unit and counit of the adjunction.
\end{remark}

\paragraph{Skeleton functors}
Given a signature $\Sigma$ of dimension $\alpha$ and an ordinal
$\beta < \alpha$, we may identify computads for $\Sigma_\beta$ with computads
for $\Sigma$ that have no generators above dimension $\beta$. More precisely, we
can define recursively a \emph{skeleton functor} and a natural transformation
\begin{align*}
  \sk_\beta^\Sigma &: \Comp_{\Sigma_\beta} \hookrightarrow \Comp_\Sigma &
  \kappa_\beta^\Sigma : \sk_\beta^\Sigma\tr_\beta^\Sigma \Rightarrow \id
\end{align*}
satisfying the following conditions
\begin{align*}
  \tr_\beta^\Sigma\sk_\beta^\Sigma
    &= \id &
  \tr_\beta^\Sigma \kappa_\beta^\Sigma
    &= \id = \kappa_\beta^\Sigma\tr_\beta^\Sigma,
\end{align*}
or equivalently that the skeleton functor is left adjoint to the truncation
functor with unit the identity and counit $\kappa_\beta^\Sigma$. In particular,
it follows that the skeleton functor is fully faithful and injective on objects.

Given a $\Sigma_\beta$-computad $C$ and an ordinal $\gamma \le \alpha$,
we may define recursively
\begin{align*}C^\gamma &=
  \begin{cases}
    \tr_\gamma^{\Sigma_\beta} C,
      &\text{if } \gamma \le \beta \\
    ((C^{\gamma'})_{\gamma'<\gamma}, (\emptyset, \{\})_{\dim i = \gamma})
      &\text{if } \beta < \gamma \le \alpha
  \end{cases}
  & \sk^{\Sigma}_\beta C &= C^\alpha,
\end{align*}
where $\{\}$ denotes the unique function out of the empty set.
Similarly, for a morphism of $\Sigma_\beta$-computads $\sigma : C\to D$, we
define recursively
\begin{align*}\sigma^\gamma &=
  \begin{cases}
    \tr_\gamma^{\Sigma_\beta} \sigma,
      &\text{if } \gamma \le \beta \\
    ((\sigma^{\gamma'})_{\gamma'<\gamma}, (\{\})_{\dim i = \gamma})
      &\text{if } \beta < \gamma \le \alpha
  \end{cases}
  & \sk^\Sigma_\beta \sigma &= \sigma^\alpha.
\end{align*}
Moreover, for a $\Sigma$\nobreakdash-computad $D$, we define the component of the counit at
$D$ recursively by
\begin{align*}\kappa^{\Sigma}_{\beta,C,\gamma} &=
  \begin{cases}
    \id_{C_\gamma},
      &\text{if } \gamma \le \beta \\
    ((\kappa^{\Sigma}_{\beta,C,\gamma'})_{\gamma'<\gamma},
    (\{\})_{\dim i = \gamma})
      &\text{if } \beta < \gamma \le \alpha
  \end{cases}
  & \kappa^{\Sigma}_{\beta,C} &= \kappa^{\Sigma}_{\beta,C,\alpha}.
\end{align*}
Naturality of the counit, as well as the two conditions above are easy to check.

\paragraph{Unbounded signatures}
We have assumed that the category $\clI$ of sorts is small. A consequence of
that is that there exists a least ordinal $\alpha$ that is greater or equal
than the dimension of every sort $i\in \clI$. The class $\Sig_\clI$ of
\emph{$\clI$-sorted signatures} is the class $\Sig_\clI(\alpha)$ of
signatures of dimension $\alpha$, and their computads and terms are
defined as above.

Observe that this definition does not really depend on
$\alpha$, since for any $\alpha' \ge \alpha$ the restriction function
$(-)^{\alpha'}_\alpha$ is a bijection and the truncation functors
$\tr_\alpha^\Sigma$ for $\Sigma\in \Sig_\clI(\alpha')$ are isomorphisms of
categories commuting with the term adjunction. Under this definition,
$\clI$-sorted signatures of arbitrary dimension $\beta$ coincide with
$\clI_\beta$-sorted signatures.

\section{Computads as presheaves}\label{sec:computads-are-presheaves}

We have defined a generalised notion of morphism of computads $\sigma : C\to D$,
where each generator of $C$ is mapped to an arbitrary term of $D$. In contrast,
the morphisms of computads usually considered are
more restricted, sending instead generators to generators~\cite{batanin_computads_1998, garner_homomorphisms_2010} . In this section, we
will identify a subcategory of computads where morphisms preserve the
generators. Generalising our  previous work~\cite{dean_computads_2022}, we will
show
that this subcategory is a presheaf topos for every signature and identify
a site of definition for it. For the rest of the section, let $\Sigma$ be an
$\clI$-sorted signature of some dimension $\alpha$.

\begin{definition}\label{def:var-to-var}
  A morphism of $\Sigma$\nobreakdash-computads $\sigma : C\to D$ is
  \emph{variable-to-variable} when $\tr_\beta^\Sigma(\sigma)$ is
  variable-to-variable for every $\beta < \alpha$, and $\sigma_i(v)$
  is a generator for
  every sort $i$ of dimension $\alpha$ and generator $v\in V_i^C$.
\end{definition}

In other words, a variable-to-variable morphism $\sigma : C\to D$ consists of
variable-to-variable morphisms $\sigma_\beta$ for all $\beta < \alpha$
satisfying the usual cocycle conditions and functions
$\sigma_i : V_i^C\to V_i^D$ for every sort $i$ of dimension $\alpha$ satisfying
the gluing condition
\[
  \phi_\delta^D\circ \sigma_i =
  \Term_{\Sigma_\beta}(\sigma_\beta)\circ \phi_\delta^C
\]
for every sort $j$ of dimension $\beta < \alpha$, and morphism
$\delta : j\to i$. variable-to-variable morphisms are closed under composition
and they contain identity morphisms, so they form a subcategory that we denote
by
\[\zeta_\Sigma : \Compv_\Sigma \to \Comp_\Sigma \]
The truncation and skeleton functors preserve the class of variable-to-variable
morphisms, so they restrict to an adjunction between the subcategories of
variable-to-variable morphisms. Moreover, for every sort $i$ of dimension at
most $\alpha$, the assignment sending a computad $C$ to the set $V_i^C$ of its
generators of sort $i$ can be extended to a functor
\[
  V_i^\bullet : \Compv_\Sigma \to \Set
\]

\begin{remark}\label{rmk:isos-are-var-to-var}
  It is easy to see that if the composition $\sigma\tau$ of two morphisms is
  variable-to-variable, then the same must be true of $\tau$. In particular,
  isomorphism of computads are variable-to-variable. Functoriality of
  $V_i^\bullet$ shows then that isomorphisms of computads induce bijections on
  the sets of generators. The converse can be shown easily by induction on the
  dimension of the signature.
\end{remark}

It is easy to see that in general the category of computads is neither
complete or cocomplete. For example, it has no terminal object, so long as the
signature contains at least one function symbol. On the contrary, the
subcategory of variable-to-variable morphisms is both complete and cocomplete,
and the inclusion functor preserves colimits and connected limits.

\begin{proposition}\label{prop:compv-cocomplete}
  The category of computads and variable-to-variable morphisms is cocomplete
  and the inclusion $\zeta_\Sigma$ is cocontinuous.
\end{proposition}

\begin{proof}
  Let $F : \clD \to \Compv_\Sigma$ a small diagram of computads and
  variable-to-variable morphisms. By induction on the dimension $\alpha$ of the
  signature,
  we may assume that colimit cocones
  $(\sigma_{\beta,d} : (Fd)_\beta\to C_\beta)_{d\in \clD}$
  have been constructed for all $\beta < \alpha$ and that they are preserved
  strictly by the truncation functors. We may then form the colimit of sets of
  generators $(\sigma_{i,d} : V^{Fd}_i \to V^C_i)_{d\in \clD}$ for every sort
  $i$ of dimension $\alpha$.

  Let $C$ the computad consisting of $C_\beta$ for every $\beta < \alpha$, the
  sets of generators $V_i^C$ for $i$ of dimension $\alpha$ and the gluing
  functions $\phi_\delta^C : V_i^C \to \Term_{\Sigma_\beta,j}(C_\beta)$ for
  every non-identity face map $\delta : j\to i$, defined by
  the universal property of the colimit on the morphisms
  \[
    V_i^{Fd}
      \xrightarrow{\phi_\delta^{Fd}}
    \Term_{\Sigma_\beta,j}((Fd)_{\beta})
      \xrightarrow{\Term(\sigma_{\beta, d})}
    \Term_{\Sigma_\beta,j}(C_{\beta}).
  \]
  where $\beta = \dim j$. The morphisms $\sigma_{\beta,d}$ and the functions
  $\sigma_{i,d}$ assemble to a cocone of variable-to-variable morphisms
  $(\sigma_d : Fd\to C)_{d\in \clD}$, whose universal property in the category
  of computads and the subcategory of variable-to-variable morphisms can be
  easily verified. Moreover, it is strictly preserved by the truncation
  functors, which concludes the induction.
\end{proof}

\begin{proposition}\label{prop:terminal-computad}
  The category of computads and variable-to-variable morphisms has a terminal
  object $\one_\Sigma$.
\end{proposition}

\begin{proof}
  Suppose that terminal computads $\one_{\Sigma_\beta}$ are given for every
  $\beta < \alpha$, strictly preserved by the truncation functors. For every
  sort $i$ of dimension $\alpha$, we may then form the limit
  \[
    V_i^{\one_\Sigma} =
    \lim_{\substack{\dim j = \beta < \alpha \\ \delta : j\to i}}
    \Term_{\Sigma_\beta,j}(\one_{\Sigma_\beta})
  \]
  The terminal computad $\one_\Sigma$ consists of the computads
  $\one_{\Sigma_\beta}$ and those sets for every sort $i$ of dimension $\alpha$.
  Its gluing functions are the obvious projections out of the limit.
\end{proof}

The existence of a terminal computad implies that the restricted functor of
terms $\Termv_\Sigma = \Term_\Sigma \circ\, \zeta_\Sigma$ can be refined to a
functor with target the slice category over $\Term_\Sigma(\one_\Sigma)$. This
slice category is equivalent to the category of presheaves on the category of
elements of $\Term_\Sigma(\one_\Sigma)$, whose objects we will call
\emph{polyplexes} following~\cite{henry_nonunital_2019}.

\begin{definition}\label{def:pplex}
  A \emph{polyplex} is a term of the terminal computad.
\end{definition}

We will denote the sort of a polyplex $p$ by $\sort(p)$. Polyplexes form a
category $\Pplex_\Sigma$, where morphisms $\delta : p\to p'$ are face maps
$\sort(p)\to \sort(p')$ satisfying that $\delta^*(p') = p$. The refinement
of the functor of terms described in the previous paragraph is given by the
functor
\begin{align*}
  \clT_\Sigma
    &: \Compv_\Sigma \times \op{\Pplex_\Sigma} \to \Set \\
  \clT_{\Sigma,p}C
    &= \{ t\in \Term_{\Sigma,\sort(p)}(C) \;:\; \Term_{\Sigma}(!)(t) = p \}
\end{align*}
where $!$ is the unique variable-to-variable morphism to the terminal computad.

\begin{proposition}\label{termv-representable}
  The functor $\clT_{\Sigma,p}$ is representable for every polyplex $p$.
\end{proposition}

\begin{proof}
  By induction on $\alpha$, we may first assume that $\clT_{\Sigma_\beta,p}$ is
  representable for every polyplex $p$ of dimension $\beta < \alpha$. If
  $\abs{p}_\beta$ is a $\Sigma_\beta$-computad representing it, then for every
  $\Sigma$\nobreakdash-computad $C$, there exists a natural isomorphism
  \[
    \clT_{\Sigma,p}C =
    \clT_{\Sigma_\beta,p}C_\beta
    \cong \Compv_{\Sigma_\beta}(\abs{p}_\beta,C_\beta)
    \cong \Compv_{\Sigma}(\sk_\beta^\Sigma\abs{p}_\beta,C),
  \]
  so $\clT_{\Sigma,p}$ is represented by $\abs{p} = \sk_\beta\abs{p}_\beta$. It
  remains to show that $\clT_{\Sigma,p}$ is representable for polyplexes $p$ of
  dimension $\alpha$ as well. We will construct such a representation
  recursively on the depth of $p$.

  Given a polyplex $p'$ of dimension less than $\alpha$, or depth less than that
  of $p$, we will denote by $\abs{p'}$ the computad representing
  $\clT_{\Sigma,p'}$. We will also denote by $t_{p'}$ the universal term in
  $\clT_{\Sigma,p'}(\abs{p'})$ inducing the representation. Finally, given a
  morphism of such polyplexes $\delta : p''\to p'$, we will denote by
  $\abs{\delta} : \abs{p''}\to\abs{p'}$ the morphism corresponding to the
  natural transformation $\clT_{\delta}$.

  Suppose first that $p$ is a polyplex of dimension $\alpha$ and depth $0$, and
  let $i$ its sort. Then there exists a family of polyplexes
  $p_\delta\in\Term_{\Sigma,j}(\one_\Sigma)$ indexed by non-identity morphisms
  $\delta : j\to i$, satisfying the usual cocycle condition, and that
  $p = \var(p_\delta)$. Then the colimit
  \[
    D = \colim_{\substack{\dim j < \alpha \\ \delta : j\to i}} \abs{p_\delta}
  \]
  is a computad with no generators of dimension $\alpha$ from the description
  of colimits in Proposition~\ref{prop:compv-cocomplete}. Let $\inc_\delta$ the
  canonical inclusion to the colimit, and let $\abs{p}$ the computad consisting
  of $D_\beta$ for all $\beta < \alpha$, has unique generator $*$ of sort $i$
  with gluing functions
  \[
    \phi^{\abs{p}}_\delta(*) =
    \clT_{\Sigma,p_\delta}(\inc_\delta)(t_{p_\delta}),
  \]
  and no other generator of dimension $\alpha$. Let also $t_{p} = \var(*)$.
  Using the universal property of the colimit defining $D$ and that $D$ has no
  top-dimensional generators, it is easy to see that evaluation at $t_p$ induces
  a natural isomorphism
  \[
    \Compv_\Sigma(\abs{p},C) \cong
    \{ v\in V_i^C \;:\; \phi_\delta^C(v)\in \clT_{\Sigma,p_\delta}C
    \text{ for all } \delta \} \cong
    \clT_pC
  \]
  for every computad $C$.

  Suppose finally that $p$ has positive depth, so that it is of the form
  $p = f[\tau]$ for some function symbol $f\in\Sigma_i$ and
  $\tau : \Cptd_\Sigma B_f\to \one_\Sigma$. We may then form the transpose
  $\tau^\dagger : B_f\to \Term_{\Sigma}(\one_\Sigma)$ of $\tau$
  under the term adjunction, which sends $b\in B_{f,j}$ to $\tau_{j}(b)$. By
  the inductive hypothesis, we may form the colimit
  \[\abs{p} = \colim_{b\in B_{f,j}} \abs{\tau^\dagger}(b)\]
  over the category of elements of $B_f$. variable-to-variable morphisms
  $\rho : \abs{p}\to C$ are in natural bijection to families of terms
  $\hat{\rho}(b)\in\clT_{\Sigma,\tau^\dagger(b)}C$ compatible with the boundary
  maps, or equivalently to morphisms $\hat{\rho} : B_f\to \Term_\Sigma(C)$
  such that $\Term(!)\hat{\rho} = \tau^\dagger$. Those
  correspond in turn to morphisms $\hat{\rho}^\dagger : \Cptd_\Sigma B_f\to C$
  such that $\tau = !\circ \hat{\rho}^\dagger$, or equivalently to terms
  $f[\hat{\rho}^\dagger]\in \clT_pC$.
\end{proof}

\begin{corollary}\label{cor:termv-fam-representable}
  The restricted functor of terms $\Termv_\Sigma$ preserves connected colimits.
\end{corollary}
\begin{proof}
  For every sort $i$, the functor $\Termv_{\Sigma,i}$ is the coproduct of the
  representable $\clT_{\Sigma,p}$ over all polyplexes $p$ of sort $i$. The
  result follows by continuity of representable functors, and commutativity of
  connected limits of sets over arbitrary coproducts.
\end{proof}

\begin{corollary}\label{cor:compv-continuous}
  The category of computads and variable-to-variable morphisms is complete and
  the inclusion $\zeta_\Sigma$ preserves connected limits.
\end{corollary}
\begin{proof}
  We have already shown that $\Compv_\Sigma$ has a terminal object, so it
  suffices to show that it has connected limits preserved by the inclusion into
  $\Comp_\Sigma$. For that, let $F : \clD\to \Compv_\Sigma$ a small, connected
  diagram of computads and variable-to-variable morphisms. By induction on the
  dimension of the signature, we may assume that limit cones
  $(\sigma_{\beta,d} : C_\beta\to (Fd)_\beta)_{d\in \clD}$ have been constructed
  for every $\beta < \alpha$ and that they are strictly preserved by the
  inclusion functors. We may then form the limits of the sets of generators
  $(\sigma_i : V_i^C \to V_i^{Fd})_{d\in \clD}$ for every sort $i$ of dimension
  $\alpha$. For every non-identity face map
  $\delta : j\to i$, we have from Corollary~\ref{cor:termv-fam-representable}
  that the morphisms
  \[
    \Term_{\Sigma_\beta,j}(\sigma_{\beta,d}) :
    \Term_{\Sigma_\beta,j}(C_\beta)\to \Term_{\Sigma_\beta,j}((Fd)_\beta),
  \]
  where $\beta = \dim j$, form a limit cone. We let
  $\phi_\delta^C : V_i^C\to \Term_{\Sigma_\beta,j}(C_\beta)$ the function
  induced by the universal property of the limit on the functions
  \[
    V_i^C \xrightarrow{\sigma_{i,d}}
    V_i^{Fd} \xrightarrow{\phi^{Fd}_\delta}
    \Term_{\Sigma_\beta,j}((Fd)_\beta).
  \]
  Let $C$ the computad consisting of $C_\beta$ for all $\beta < \alpha$, the
  set $V_i^C$ for every sort $i$ of dimension $\alpha$ and the gluing functions
  above. The morphisms $\sigma_{\beta,d}$ and the functions $\sigma_{i,d}$
  give rise to a cone $(\sigma : C\to Fd)_{d\in \clD}$, whose universal property
  in the subcategory of variable-to-variable morphisms can be verified
  immediately. Its universal property in the category of computads and all
  morphisms follows then easily by preservation of connected limits by
  $\Termv_\Sigma$, Finally, this limit cone is strictly preserved by the
  truncation functors, which concludes the induction.
\end{proof}

At this point, we have all the ingredients needed to show that $\Compv_\Sigma$
is a presheaf topos. We define a \emph{plex} to be a generator of the terminal
computad $\one_\Sigma$. Plexes form a direct subcategory $\Plex_\Sigma$ of
$\Compv_\Sigma$, where
\[
  \Plex_\Sigma(p,p')
  = \Compv_\Sigma(\abs{p},\abs{p'})
  \cong \clT_{\Sigma,p}(\abs{p'})
\]
and $\dim p = \dim(\sort(p))$. Plexes familially represent the functors sending
a computads to its generators, in the sense that there exist natural
isomorphisms
\[
  V_i^C \cong
  \coprod_{\substack{p\in \Plex_\Sigma \\ \sort(p) = i}}
  \Compv_\Sigma(\abs{p}, C).
\]
From the construction of colimits of variable-to-variable morphisms in
Proposition~\ref{prop:compv-cocomplete}, it easy to see that the functors
represented by the plexes are cocontinuous. Moreover, they jointly reflect
isomorphisms by Remark~\ref{rmk:isos-are-var-to-var}. The following theorem is then an immediate
consequence of~\cite[Proposition~5.14]{dean_computads_2022}.

\begin{theorem}\label{thm:compv-topos}
  The nerve functor $N : \Compv_\Sigma\to [\op{\Plex_\Sigma}, \Set]$ defined by
  \[
    (NC)(p) =
    \Compv_\Sigma(\abs{p},C)
    \cong \{ v\in V^C_{\sort(p)} \;:\; \Term_\Sigma(!)(\var v) = p \},
  \]
  is an equivalence of categories.
\end{theorem}

\section{Algebras over a signature}\label{sec:algebras}

In this section, we introduce the semantics of $\clI$-sorted signatures.
We define algebras for $\Sigma$ to be algebras for the term monad $\M_\Sigma$, and
give a simpler description of them in terms of presheaves equipped with a
function for every function symbol of $\Sigma$, satisfying certain boundary
axioms dictated by the boundary terms. By definition, computads gives rise to
free algebras, and we will show that our generalised morphisms of
computads are precisely the morphisms between the algebras they
generate. As before, in this section, $\Sigma$ denotes an $\clI$-sorted
signature of some dimension $\alpha$.

\begin{definition}\label{def:algebra}
  A \emph{$\Sigma$\nobreakdash-algebra} is an algebra for the monad
  $(\M_\Sigma, \eta_\Sigma, \mu_\Sigma)$ induced by the term adjunction
  $\Cptd_\Sigma\dashv\Term_\Sigma$. We will denote their category by
  $\Alg_\Sigma = \Alg_{\M_\Sigma}$.
\end{definition}

An algebra $\bbX = (X, u_\bbX)$ consists therefore of a \emph{carrier}
presheaf $X$ on $\clI_\alpha$, and a morphism
$u^\bbX : \M_\Sigma(X)\to X$, the \emph{$\M_\Sigma$-action}, satisfying the
following unit and associativity axioms
\begin{align*}
  u^\bbX \circ\eta_\Sigma &= \id &
  u^\bbX\circ\mu_\Sigma &= u^\bbX\circ \M_\Sigma (u^\bbX).
\end{align*}
The unit axiom prescribes the value of the action on generators by
\[ u^\bbX(\var x) = x \]
for every $i\in \clI_\alpha$ and $x\in X_i$. The associativity axiom, on the
other hand, gives rise to a recursive formula for composite terms: Given a
composite term $t = f[\tau]$, we may form the transpose
$\tau^\dagger : B_f\to \M_\Sigma(X)$ of $\tau$ under the term adjunction. By
definition, $\tau = \varepsilon_{\Sigma,X}\circ\Cptd_\Sigma(\tau^\dagger)$, so
\begin{equation}\label{eq:recursive-action}
  u^\bbX(f[\tau]) =
  u^\bbX\mu_\Sigma(f[\Cptd_\Sigma(\tau^\dagger)]) =
  u^\bbX(f[\Cptd_\Sigma(u^\bbX\tau^\dagger)]).
\end{equation}
We define the \emph{interpretation} of a function symbol $f\in \Sigma_i$ of
dimension $\beta\le\alpha$ in $\bbX$ to be the function
\[f^\bbX : [\op{\clI_\beta},\Set](B_f,\tr_\beta X) \to X_i \]
given on a morphism $\sigma : B_f\to \tr_\beta X$ by
\[f^\bbX(\sigma) = u^\bbX(f[\Cptd_\Sigma(\sigma)]).\]
A recursive formula shows that an algebra $\bbX$ is uniquely determined by its
carrier presheaf and the interpretations of the function symbols. Moreover,
those interpretation functions can be freely chosen, so long as certain boundary
condition is satisfied.

\begin{proposition}\label{prop:interpret-determine-obj}
  Let $X$ a presheaf on $\clI_\alpha$. Actions $u : \M_\Sigma(X)\to X$ are in
  bijection to families of
  \begin{itemize}
    \item actions
      $u_\beta : \M_{\Sigma_\beta}(\tr_\beta X)\to \tr_\beta X$
      for $\beta < \alpha$
    \item functions $\hat{f} : [\op{\clI},\Set](B_f,X) \to X_i$ for $i$
      of dimension $\alpha$ and $f\in \Sigma_i$,
  \end{itemize}
  satisfying the usual cocycle conditions, and the following boundary condition
  \[ \delta^*\hat{f} = u_\beta(\M_{\Sigma_\beta}(-)(t_{f,\delta})) \]
  for non-identity face maps $\delta : j\to i$, where $\beta = \dim j$.
\end{proposition}
\begin{proof}
  As explained above, the morphisms $u_\beta$ and the functions $\hat{f}$
  determine $u$ uniquely. Conversely, given such morphisms and functions, we
  can build such $u$ recursively. We first define $u$ on terms of dimension
  $\beta < \alpha$ to coincide with $u_\beta$ and define it on generators
  $x\in X_i$ of dimension $\alpha$ by
  \[u(\var x) = x.\]
  This assignment is compatible with the boundary maps by the cocycle conditions
  and the fact that each $u_\beta$ satisfies the unit axiom. We define then $u$
  on composite term $t = f[\tau]$ of dimension $\alpha$ by the recursive
  formula~(\ref{eq:recursive-action}):
  \[u(f[\tau]) = \hat{f}(u\circ \tau^\dagger). \]
  This assignment is compatible with face maps by the boundary condition,
  hence it defines a morphism $u : \M_\Sigma (X)\to X$.

  This morphism satisfies the associativity axiom on terms of dimension less
  than $\alpha$, and the unit axiom on all terms. The unit axiom also implies
  the associativity axiom for generators of dimension $\alpha$, so let
  $t = f[\rho]$ a composite term of $\Cptd_\Sigma\M_\Sigma(X)$ of dimension
  $\alpha$. We may assume that the associativity axiom
  holds in the image of the transpose $\rho^\dagger : B_f\to
  \M_\Sigma\M_\Sigma X$ of $\rho$ by induction on depth, and compute that
  \begin{align*}
    u\circ \M_\Sigma(u)(t)
      &= u(f[(\Cptd_\Sigma u) \circ \rho]) \\
      &= \hat{f}(u\circ ((\Cptd_\Sigma u) \circ \rho)^\dagger) \\
      &= \hat{f}(u\circ \M_\Sigma(u) \circ \rho^\dagger) \\
      &= \hat{f}(u\circ \mu_\Sigma \circ \rho^\dagger) \\
      &= \hat{f}(u\circ (\varepsilon_\Sigma\rho)^\dagger) \\
      &= u(f[\varepsilon_\Sigma\rho]) \\
      &= u\circ \mu_\Sigma(t)
  \end{align*}
  Therefore, the associativity axiom holds for all terms, and $\bbX = (X,u)$ is
  an algebra.

  It remains to show that this algebra gives rise to the data that we started
  from. The cocycle conditions imply that $\tr_\beta u = u$ for every
  $\beta < \alpha$. Let therefore $f\in \Sigma_i$ of dimension $\alpha$ and
  $\tau : B_f\to X$. By definition of transposition and the unit axiom, we have
  that
  \[
    u\circ \Cptd_\Sigma(\tau)^\dagger =
    u\circ \M_\Sigma(\tau)\circ \eta_\Sigma =
    u\circ \eta_\Sigma\circ \tau =
    \tau
  \]
  and hence that
  \[
    f^\bbX(\tau) =
    u(f[\Cptd_\Sigma(\tau)]) =
    \hat{f}(u\circ\Cptd_\Sigma(\tau)^\dagger) =
    \hat{f}(\tau).
  \]
  We see that the interpretation of $f$ in this algebra coincides with the
  function $\hat{f}$, proving the bijection of the proposition.
\end{proof}

Morphisms of algebras $f : \bbX\to \bbY$ are the morphisms
$\sigma : X\to Y$ satisfying that
\[
  u^\bbY\circ (\M_\Sigma \sigma) = \sigma \circ u^\bbX.
\]
This equation easily implies that morphisms of algebras must preserve the
interpretation functions
\begin{equation}\label{eq:compatibility-interpretations}
  \sigma\circ f^\bbX = f^\bbY(\tr_\beta\sigma\circ -)
\end{equation}
for every $f\in \Sigma_i$ of dimension $\beta \le \alpha$. The recursive
formula~(\ref{eq:recursive-action}) immediately implies that the converse also holds.

\begin{proposition}\label{prop:interpret-determine-mor}
  A morphism of algebras $\sigma : \bbX\to \bbY$ is a morphism between their
  underlying presheaves $\sigma : X\to Y$ satisfying
  equation~(\ref{eq:compatibility-interpretations}) for every function symbol of
  $\Sigma$.
\end{proposition}

In the rest of the section, we will study the connection between
computads and algebras. The category of algebras is equipped by definition with
an adjunction $\F_\Sigma\dashv \U_\Sigma$ to the
category of presheaves on $\clI_\alpha$ and with a functor $\K_\Sigma$ from
the category of computads, making the following triangle commute in
both directions.
\[
  \begin{tikzcd}
    & {[\op{\clI_\alpha},\Set]} \\
    \Comp_\Sigma && \Alg_\Sigma
    \arrow[""{name=0, anchor=center, inner sep=0}, "{\Cptd_\Sigma}",
    shift left=1, from=1-2, to=2-1]
    \arrow[""{name=1, anchor=center, inner sep=0}, "{\Term_\Sigma}",
    shift left=2, from=2-1, to=1-2]
    \arrow[""{name=2, anchor=center, inner sep=0}, "{\F_\Sigma}"',
    shift right=1, from=1-2, to=2-3]
    \arrow[""{name=3, anchor=center, inner sep=0}, "{\U_\Sigma}"',
    shift right=2, from=2-3, to=1-2]
    \arrow["{\K_\Sigma}"', from=2-1, to=2-3]
    \arrow["\dashv"{anchor=center, rotate=55}, draw=none, from=2, to=3]
    \arrow["\dashv"{anchor=center, rotate=123}, draw=none, from=0, to=1]
  \end{tikzcd}
\]
The forgetful functor $\U_\Sigma$ sends an algebra to its carrier presheaf,
while the comparison functor $\K_\Sigma$ sends a computad $C$ to its presheaf
of terms $\Term_\Sigma(C)$ equipped with the action
$u^C=\Term_\Sigma(\varepsilon_{\Sigma,C})$. Henceforth,
we will treat $\K_\Sigma$ as a subcategory inclusion and suppress it in the
notation. This abuse of notation partially justified by Proposition~\ref{prop:comparison-fully-faithful}.

\begin{lemma}\label{lem:generators-determine}
  If morphisms $\sigma, \tau : C\to \bbX$ from a computad to
  an algebra agree on all generators of $C$, then they are equal.
\end{lemma}

\begin{proof}
  By induction on the dimension $\alpha$ of the signature, we may assume that
  $\tr_\beta\sigma = \tr_\beta\tau$ for all $\beta < \alpha$, so it remains to
  show that they agree on composite terms of dimension $\alpha$ as well. Given
  a term $t = f[\rho]$, we may further assume by structural induction that the
  morphisms agree on terms in the image of the transpose $\rho^\dagger$ as well.
  Using that, we see that
  \begin{align*}
    \sigma(f[\rho])
      &= \sigma\circ u^{C}
        (f[\Cptd_\Sigma\rho^\dagger])\\
      &= u^\bbX \circ \M_\Sigma(\sigma)(f[\Cptd_\Sigma\rho^\dagger]) \\
      &= u^\bbX (f[\Cptd_\Sigma(\sigma\rho^\dagger)]) \\
      &= u^\bbX (f[\Cptd_\Sigma(\sigma'\rho^\dagger)]) \\
      &= \dots = \tau(f[\rho]),
 \end{align*}
  so the two morphisms are equal.
\end{proof}

\begin{proposition}\label{prop:comparison-fully-faithful}
  The functor $\K_\Sigma$ is fully faithful.
\end{proposition}
\begin{proof}
  The composite $\Term_\Sigma = \U_\Sigma \K_\Sigma$ is clearly faithful, so the
  same must hold of $\K_\Sigma$. To show that it is full, let $C, D$ be
  computads and $\sigma: \K_\Sigma C\to \K_\Sigma D$ a morphism of
  algebras. By induction on dimension, we may assume that there
  exists a unique morphisms of computads
  $\rho_\beta : C_\beta\to D_\beta$ for every $\beta < \alpha$, such that
  \begin{equation}\label{eq:comparison-ff-equation}
    \K_{\Sigma_\beta}(\rho_\beta) = \tr_\beta\sigma.
  \end{equation}
  Uniqueness of those morphisms implies that the usual cocycle conditions are
  satisfied. Let $\rho : C\to D$ consist of $\rho_\beta$ for all $\beta<\alpha$
  and the functions
  \[\rho_i = \sigma_i \circ \var : V_i^C\to \Term_{\Sigma,i}(D)\]
  for every sort $i$ of dimension $\alpha$. Using equation~(\ref{eq:comparison-ff-equation}), it is easy to see that $\rho$ is a
  well-defined morphism of computads. The morphisms of algebras
  $\sigma$ and $\K_\Sigma(\rho)$ agree on generators, so they must be equal by
  Lemma~\ref{lem:generators-determine}.
\end{proof}

\begin{proposition}\label{prop:universal-property-computads}
  Morphisms $\sigma : C\to \bbX$ from a computad to an algebra $\bbX$ are in
  bijection to
  \begin{itemize}
    \item morphisms $\rho_\beta : \tr_\beta C\to \tr_\beta \bbX$ for every
      $\beta < \alpha$,
    \item functions $\rho_i : V_i^C\to X_i$ for every sort $i$ of dimension
      $\alpha$,
  \end{itemize}
  satisfying the usual cocycle conditions, and the boundary condition
  \[\delta^*\rho_i = \rho_{\beta,j}\phi^C_\delta\]
  for non-identity face maps $\delta : j \to i$, where $\beta = \dim j$.
\end{proposition}

\begin{proof}
  A morphism $\sigma$ gives rise to such data by letting
  $\rho_\beta =\tr_\beta\sigma$ for $\beta < \alpha$ and letting
  $\rho_i = \sigma_i\circ \var$ for $i$ of dimension $\alpha$. Moreover, this
  assignment is injective by Lemma~\ref{lem:generators-determine}, so it
  remains to construct a morphism $\sigma : C\to \bbX$ out of that data.

  On terms of dimension $\beta < \alpha$, we let $\sigma$ coincide with
  $\rho_\beta$, while on generators of sort $i$ of dimension $\alpha$, we let
  $\sigma$ coincide with $\rho_i$. this assignment is compatible with the
  face maps by the cocycle and boundary conditions. Given a composite term
  $t = f[\tau]$ of sort $i$, we may assume recursively on depth that $\sigma$
  has been defined on the image of the transpose $\tau^\dagger$ and let
  \[\sigma(f[\tau]) = f^\bbX(\sigma\tau^\dagger).\]
  Given a non-identity face map $\delta : j\to i$, we let $\beta = \dim j$ and
  calculate that
  \begin{align*}
    \delta^*\sigma(f[\tau])
      &= \delta^*f^\bbX(\sigma\tau^\dagger) \\
      &= u^\bbX\circ (\M_\sigma\sigma)\circ
        (\M_\Sigma\tau^\dagger)(t_{f,\delta}) \\
      &= \sigma\circ \Term_\Sigma(\varepsilon_{\Sigma,C})\circ
        (\M_\Sigma\tau^\dagger)(t_{f,\delta}) \\
      &= \sigma \circ \Term_\Sigma(\tau)(t_{f,\delta}) \\
      &= \sigma(\delta^*(f[\tau])),
  \end{align*}
  so this assignment defines a morphism of presheaves
  $\sigma : \Term_\Sigma C\to X$, which restricts to the morphisms $\rho_\beta$
  and the functions $\rho_i$. It remains to show that $\sigma$ is a morphism of
  algebras. By Proposition~\ref{prop:interpret-determine-mor}, it
  suffices to show that it preserves the interpretations of every function
  symbol of $\Sigma$. This is the case for function symbols of dimension
  $\beta < \alpha$, since $\rho_\beta$ is a morphism of
  algebras. On the other hand, let $f\in \Sigma_i$ of
  dimension $\alpha$ and $\hat{\tau} : B_f\to \Term_\Sigma(C)$ a morphism. Then
  \begin{align*}
    \sigma(f^C(\hat{\tau}))
      &= \sigma(f[\varepsilon_{\Sigma,C}\Cptd_\Sigma(\hat{\tau})]) \\
      &= f^\bbX(\sigma(\varepsilon_{\Sigma,C}
        \Cptd_\Sigma(\hat{\tau}))^\dagger)\\
      &= f^\bbX(\sigma\hat{\tau}^{\dagger\dagger}) = f^\bbX(\sigma\hat{\tau})
  \end{align*}
  so $\sigma$ is a morphism of algebras.
\end{proof}

\begin{corollary}\label{cor:universal-property-truncated-computads}
  Morphisms $\sigma : C\to \bbX$ from a computad with no generators of
  dimension at least $\beta$ to an arbitrary algebra $\bbX$ are in
  bijection to families of morphisms $\sigma_\gamma : \tr_\gamma^\Sigma C\to
  \tr_\gamma^\Sigma\bbX$ for every $\gamma < \beta$ satisfying the usual cocycle
  conditions.
\end{corollary}

\section{Cofibrancy of computads}\label{sec:computads-are-cofibrant}

Computads have recently proven useful in the study of the homotopy theory of
$\omega$\nobreakdash-categories. The inclusion of the free $\omega$\nobreakdash-category on a sphere
into the free $\omega$\nobreakdash-category on a disk cofibrantly generate a weak
factorisation system in the category of $\omega$\nobreakdash-categories and strict morphisms
for which computads are cofibrant. Moreover, there exists an adjunction between
computads and $\omega$\nobreakdash-categories~\cite{batanin_computads_1998} generating the
\emph{universal cofibrant replacement} comonad for this factorisation system~\cite{garner_homomorphisms_2010}. For strict $\omega$\nobreakdash-categories,
this weak factorisation system is part of a model structure~\cite{lafont_folk_2010}, and it has been shown that conversely every cofibrant
strict $\omega$\nobreakdash-category is free on a computad~\cite{metayer_cofibrant_2008}.

Most of those results hold for our computads for $\clI$-sorted theories
verbatim. Proposition~\ref{prop:universal-property-computads} describes a
universal property of algebras free on a computad, analogous to the universal
property of free $\omega$\nobreakdash-categories~\cite{street_algebra_1987}, which can
be used to construct a right adjoint to the free algebra functor
\[
  \Free_\Sigma :
  \Compv_\Sigma\xhookrightarrow{\zeta_\Sigma}
  \Comp_\Sigma \xhookrightarrow{\K_\Sigma}
  \Alg_\Sigma.
\]
We will show that the comonad induced by this adjunction is the universal
cofibrant replacement for certain weak factorisation system, and that algebras
free on a computad are the cofibrant objects.

To set the notation, fix an $\clI$-sorted signature $\Sigma$ of some dimension
$\alpha$. The \emph{representable computad} on a sort $i\in \clI_\alpha$ is the
computad
\[ \bbD^i_\Sigma = \Cptd_\Sigma(\clI_\alpha(-,i)).\]
Its \emph{boundary} is the computad obtained by removing its top-dimensional
generators, or equivalently by
\[
  \partial\bbD^i_\Sigma =
  \Cptd_\Sigma\left(\colim_{\delta : j\to i}
    \clI_\alpha(-,j)\right)
\]
where the colimit is over all non-identity morphisms $\delta : j\to i$. The
\emph{boundary inclusion}
\[
  \iota_{\Sigma,i} : \partial\bbD^i_\Sigma \hookrightarrow \bbD_\Sigma^i
\]
is the morphism of computads induced by the morphism of presheaves whose
component at $\delta$ is $\delta_* = (\delta\circ -)$. The following proposition
is an immediate consequence of the Yoneda lemma and the term adjunction.

\begin{proposition}\label{prop:term-representable}
  The computad $\bbD^i$ represents the functor $\Term_{\Sigma,i}$
  of terms of sort $i$. Its boundary $\partial \bbD^i_\Sigma$
  represents the functor of types of sort $i$
  \[ \Type_{\Sigma,i} = \lim_{\delta : j\to i} \Term_{\Sigma,j} \]
  and the boundary inclusion $\iota_{\Sigma,i}$ induces the natural
  transformation sending a term $t$ to the family of terms $(\delta^*t)$.
\end{proposition}

Using instead the adjunction $\F_\Sigma\dashv \U_\Sigma$, it is easy to see that
the functor sending an algebra $\bbX$ to the set $X_i$ is represented
by the free algebra on $\bbD^i_\Sigma$, and that the functor
represented by the free algebra on $\partial\bbD^i_\Sigma$ admits a
similar description.

The set of boundary inclusions $I = \{\iota_{\Sigma,i} \;:\; i\in \clI_\alpha\}$
cofibrantly generates a weak factorisation system $(\clL,\clR)$ in
the category of algebras via the small object argument, in light of
Corollary~\ref{cor:algebras-presentable}. We will call morphisms of $\clL$
\emph{cofibrations} and morphism of $\clR$ \emph{trivial fibrations}. We will
call an algebra $\bbX$ \emph{cofibrant} when the unique morphism
$\emptyset\to \bbX$ from the initial algebra is a cofibration. We
will also say that a morphism $\bbX\to \bbY$ is a \emph{cofibrant replacement}
when it is a trivial fibration and $\bbX$ is cofibrant.

\begin{proposition}\label{prop:free-implies-cofibrant}
  Free algebras on a computad are cofibrant.
\end{proposition}

\begin{proof}
  We will show that the unique morphism $\emptyset \to C$ is a transfinite
  composition of pushouts of coproducts of boundary inclusions in the category
  of algebras for every computad $C$. Closure of cofibrations under
  those operations is a standard result~\cite{hovey_model_2007}. In order to show that, we will introduce a variant
  of the skeleton functors.

  Let $\sk_{\partial \beta}^\Sigma C$ the computad obtained by removing all
  generators of $C$ of dimension at least $\beta$ for every $\beta \le \alpha+1$
  \[
    \sk_{\partial \beta}^\Sigma C =
    \sk_{\beta}^\Sigma ( (\tr_\gamma^\Sigma C)_{\gamma < \beta},
      (\emptyset, \{\})_{\dim j = \beta}).
  \]
  Inclusions of the sets of generators give rise to variable-to-variable
  morphisms
  \[
    \kappa_\gamma^\beta :
    \sk_{\partial \gamma}^\Sigma C \to
    \sk_{\partial \beta}^\Sigma
  \]
  for every $\gamma \le \beta$ satisfying the obvious cocycle conditions. That
  data defines a chain of variable-to-variable morphisms which is cocontinuous
  by Proposition~\ref{prop:compv-cocomplete} and has transfinite composite
  $\emptyset \to C$. It remains to show that this chain remains cocontinuous
  when viewed as a chain of algebras and that the morphisms
  $\kappa_\beta^{\beta+1}$ are pushouts of coproducts of boundary inclusions
  for every $\beta \le \alpha$.

  For the former, let $\beta\le\alpha$ a limit ordinal and let
  \[\tau^\gamma : \sk_{\partial \gamma}^\Sigma C\to \bbX\]
  for $\gamma < \alpha$ a cocone under the restriction of this chain to $\beta$.
  Then the morphisms
  \[
    \tau_\gamma = \tr_\gamma^\Sigma(\tau^{\gamma+1}):
    \tr_\gamma^\Sigma\sk_{\partial \beta}^\Sigma C\to \bbX
  \]
  satisfy the cocycle conditions and give rise to a morphism
  $\tau : \sk_{\partial \beta}^\Sigma C\to \bbX$ by Corollary~\ref{cor:universal-property-truncated-computads}.
  By Lemma~\ref{lem:generators-determine} it follows immediately that $\tau$ is the
  unique morphism such that $\tau^\gamma = \tau\kappa_\gamma^\beta$ for every
  $\gamma < \beta$. Therefore, this chain is cocontinuous in when viewed as a
  chain of algebras.

  For the latter, given an ordinal $\beta \le \alpha$, we may form the
  commutative square
  \[\begin{tikzcd}
    \coprod\limits_{\substack{\dim i = \alpha \\ v\in V_i^C}}
      \partial \bbD^i_\Sigma &
    \coprod\limits_{\substack{\dim i = \alpha \\ v\in V_i^C}}\bbD^i_\Sigma \\
    {\sk^\Sigma_{\partial\beta}C} &
    {\sk^\Sigma_{\partial(\beta+1)}C}
    \arrow["(\iota_{\Sigma,i})", from=1-1, to=1-2]
    \arrow["\phi"', from=1-1, to=2-1]
    \arrow["{\kappa^{\beta+1}_\beta}"', from=2-1, to=2-2]
    \arrow["\psi", from=1-2, to=2-2]
  \end{tikzcd}\]
  where $\psi$ classifies the generators of $C$ of dimension $\beta$ under the
  bijection of Proposition~\ref{prop:term-representable}, and $\phi$ classifies
  their boundary types. The left adjoints $\F_\Sigma$ preserves colimits, while
  $\Cptd_\Sigma$ reflects them, so both coproducts of the square are also
  coproducts of algebras. Using Corollary
 ~\ref{cor:universal-property-truncated-computads} for the computads on the
  bottom row and Proposition~\ref{prop:universal-property-computads}, we see
  that this square is a pushout of algebras.
\end{proof}

In the proof of the proposition above, we described each computad as a
transfinite composite of pushouts of coproducts of boundary inclusions. This
descriptions allows us to define a right adjoint
\[
  \Und_\Sigma : \Alg_\Sigma \to \Compv_\Sigma
\]
to the free functor $\Free_\Sigma$ sending a computad to the
algebra it generates. This adjunctions defines a comonad on the
category of algebras, whose underlying pointed endofunctor we will
denote by
\begin{align*}
  \Cof_\Sigma &: \Alg_\Sigma \to \Alg_\Sigma \\
  r_\Sigma &: \Cof_\Sigma \Rightarrow \id.
\end{align*}
The definition of $\Und_\Sigma$ is recursive on the dimension $\alpha$ of the
signature, so suppose that the adjunctions $\Free_{\Sigma_\beta}\dashv
\Und_{\Sigma_\beta}$ have been defined for all $\beta < \alpha$ and that they
are compatible with the truncation functors.

The underlying computad $\Und_\Sigma \bbX$ of an algebra $\bbX$
consists of the computads $\Und_{\Sigma_\beta}\tr_\beta^\Sigma\bbX$ for
$\beta < \alpha$ and the sets of generators and gluing functions defined by the
following pullback square
\[
  \begin{tikzcd}[row sep = large]
	{V_i^{\Und_\Sigma\bbX}} && {X_i} \\
	{\lim\limits_{\substack{\dim j = \beta < \alpha \\\delta : j\to i}}
  (\Cof_{\Sigma_{\beta}}\tr^\Sigma_{\beta}\bbX)_j} &
  {\lim\limits_{\substack{\dim j = \beta < \alpha \\\delta : j\to i}}
  (\tr^\Sigma_{\beta}\bbX)_j} &
  {\lim\limits_{\substack{\dim j = \beta < \alpha \\\delta : j\to i}} X_j}.
	\arrow[equals, from=2-2, to=2-3]
	\arrow["{(\phi_{\delta}^{\Und_\Sigma\bbX})}"', dashed, from=1-1, to=2-1]
	\arrow[dashed, from=1-1, to=1-3]
	\arrow["{(\delta^*)}", from=1-3, to=2-3]
	\arrow["{(r_{\Sigma_{\beta}})}"', from=2-1, to=2-2]
\end{tikzcd}\]
for every sort $i$ of dimension $\alpha$. The variable-to-variable morphism
$\Und_\Sigma\sigma$ induced by a morphism $\sigma : \bbX\to \bbY$ consists
similarly of the morphisms $\Und_{\Sigma_\beta}\tr_\beta^\Sigma\sigma$ for
$\beta < \alpha$ and the functions
$V_i^{\Und_\Sigma\bbX}\to V_i^{\Und_\Sigma\bbY}$ for $i$ of dimension $\alpha$
with components
\[
  (\Und_\Sigma\sigma)_i((t_\delta), x) =
  (((\Cof_{\Sigma_{\beta}}\tr^\Sigma_{\beta}\sigma)(t_\delta)),\sigma(x)).
\]
This assignment is functorial and commutes clearly with the truncation functors.
The counit of the adjunction $r_{\Sigma,\bbX} : \Cof_\Sigma\bbX\to \bbX$ is
the morphism corresponding to the morphisms
$r_{\Sigma_\beta,\tr_\beta^\Sigma\bbX}$ and the projection functions
$V_i^{\Und_\Sigma \bbX} \to X_i$ under the bijection of Proposition~\ref{prop:universal-property-computads}.

\begin{proposition}\label{prop:free-und-adjunction}
  The functor $\Und_\Sigma$ is right adjoint to $\Free_\Sigma$.
\end{proposition}

\begin{proof}
  Variable-to-variable morphisms $\sigma : C\to \Und_\Sigma\bbX$ consist of
  variable-to-variables morphisms
  $\sigma_\beta : C_\beta\to \Und_{\Sigma_\beta}\tr_\beta^\Sigma\bbX$ satisfying
  the usual cocycle conditions and a pair of functions
  \begin{align*}
    \sigma_{i,1}
      &: V_i^C\to \lim_{\substack{\dim j = \beta < \alpha \\\delta : j\to i}}
        (\tr_{\Sigma_{\beta}}\bbX)_j \\
    \sigma_{i,2}
      &: V_i^C\to X_i
  \end{align*}
  satisfying gluing conditions and that $(\sigma_{i,1},\sigma_{i,2})$ defines a
  function into the pullback $V_i^{\Und_\Sigma \bbX}$. The gluing condition is
  equivalent to
  \[
    \sigma_{i,1} =
    (\Term_{\Sigma_{\beta}}(\sigma_{\beta})
    \circ\phi_\delta^C)_{\delta : j\to i},
  \]
  in the presence of which the other condition becomes
  \[
    \delta^*\sigma_{i,2} =
    r_{\Sigma_{\beta}}\Term_{\Sigma_{\beta},j}(\sigma_{\beta})\phi_{\delta}^C
  \]
  for non-identity face map $\delta : j\to i$ where $\beta = \dim j$. By the
  inductive hypothesis, the morphisms $\sigma_\beta$
  are in bijection to morphisms
  $\sigma_\beta^\dagger : C_\beta\to \tr_\beta^\Sigma\bbX$ satisfying the same
  cocycle conditions. Under this bijection, the condition above becomes
  \[
    \delta^*\sigma_{i,2} =
    \sigma_{\beta,j}^\dagger\phi_{\delta}^C,
  \]
  so the morphisms $\sigma_\beta^\dagger$ and the functions $\sigma_{i,2}$
  determine uniquely a morphism $\sigma^\dagger : C\to \bbX$ by Proposition~\ref{prop:universal-property-computads}.
  Using Lemma~\ref{lem:generators-determine}, it is easy to see that this bijection is given
  by $\sigma^\dagger = r_{\Sigma,\bbX}\circ \Free_\Sigma\sigma$. Therefore, the
  bijection is natural and $r_\Sigma$ is the counit of the adjunction.
\end{proof}

Recall that trivial fibrations are the morphisms in the right class of the
factorisation system generated by the boundary inclusions, so they are morphisms
of algebras $\sigma : \bbX\to \bbY$ such that every commutative
square of the form
\[\begin{tikzcd}
  {\partial\bbD^i_\Sigma} & \bbX \\
  {\bbD^i_\Sigma} & \bbY
  \arrow["T", from=1-1, to=1-2]
  \arrow["t", from=2-1, to=2-2]
  \arrow[hook, from=1-1, to=2-1]
  \arrow["{\sigma}", from=1-2, to=2-2]
  \arrow["{v(i,T,t)}"{description}, dashed, from=2-1, to=1-2]
\end{tikzcd}\]
admits a diagonal lift $v(i,T,t)$. An \emph{algebraic trivial fibration} is
instead a trivial fibration with a choice of lifts for every such square.
Morphisms of algebraic trivial fibrations with common target $\bbY$ are
morphisms in the slice category $\Alg_\Sigma/\bbY$ preserving the lifts.

The counit $r_{\Sigma,\bbX} : \Cof_\Sigma \bbX \to \bbX$ can be equipped with
the structure of an algebraic trivial fibration as follows: A commutative
square of the form
\[\begin{tikzcd}
  {\partial\bbD^i_\Sigma} & \Cof_\Sigma\bbX \\
  {\bbD^i_\Sigma} & \bbX
  \arrow["T", from=1-1, to=1-2]
  \arrow["t", from=2-1, to=2-2]
  \arrow[hook, from=1-1, to=2-1]
  \arrow["{r_{\Sigma,\bbX}}", from=1-2, to=2-2]
\end{tikzcd}\]
amounts to the choice of compatible terms $T = (t_\delta)_{\delta : j\to i}$ of
$\Cof_{\Sigma}\bbX$ together with an element $t\in X_i$ satisfying the
compatibility condition
\[r_{\Sigma,\bbX}(t_\delta) = \delta^*t,\]
or equivalently a generator $(T,t)$ of $\Cof_\Sigma\bbX$ of sort $i$. A lift for
this square is then given by the morphism corresponding to the term
\[v(i,T,t) = \var(T,t).\]
It is not hard to see that $(r_{\Sigma,\bbX},v)$ is initial among algebraic
trivial fibrations with target $\bbX$. Given any algebraic trivial fibration
$(\sigma : \bbY\to \bbX, v')$, we may define a morphism
$\tau : \Cof_\Sigma\bbX\to \bbY$ of algebraic trivial fibrations recursively
by letting for every $\beta \le  \alpha$, the morphism
$\tau_\beta : \tr_\beta^\Sigma\Cof_\Sigma\bbX\to \tr_\beta^\Sigma\bbY$
correspond under the bijection of Proposition~\ref{prop:universal-property-computads} to the morphisms $\tau_\gamma$ for
$\gamma < \beta$ and the functions $\tau_i : V_i^{\Und_\Sigma\bbX}\to Y_i$
given by
\[
  \tau_i((t_\delta)_{\delta : j\to i},t) =
  v'(i,(\tau_{\beta}(t_\delta))_{\delta : j\to i},t).
\]
The morphism $\tau = \tau_\alpha : \Cof_{\Sigma}\bbX\to \bbY$ is a morphism of
algebraic trivial fibrations, since $r_{\Sigma,\bbX} = \sigma\tau$ and
\[\tau \circ v(i,T,t) = v'(i,\tau\circ T, t)\]
for every commutative square as above. Moreover, $\tau$ is unique by Lemma~\ref{lem:generators-determine}. This observation combined with the recognition
criterion~\cite[Proposition~2.6]{garner_homomorphisms_2010} show the following
corollary.

\begin{corollary}\label{cor:universal-cofibrant-replacement}
  The pointed endofunctor $(\Cof_\Sigma, r_\Sigma)$ underlies the universal
  cofibrant comonad for the weak factorisation system cofibrantly generated by
  the boundary inclusions.
\end{corollary}

The existence of a cofibrant replacement functor taking values in free
algebras together with Cauchy completeness of the category of
computads, shown in Corollary~\ref{cor:comp-cauchy} allow us to
prove the converse of Proposition~\ref{prop:free-implies-cofibrant} by the
same argument used for strict $\omega$\nobreakdash-categories~\cite{metayer_cofibrant_2008}.

\begin{corollary}\label{cor:cofibrant-implies-free}
  Cofibrant algebras are free on a computad.
\end{corollary}
\begin{proof}
  Let $\bbX$ a cofibrant algebra. By Corollary~\ref{cor:universal-cofibrant-replacement}, there exists a computad $C$ and a
  trivial fibration $r: C\to \bbX$. Since $\bbX$ is cofibrant, the map $r$
  admits a section $s : \bbX\to C$. Being a section, the endomorphism
  $sr : C\to C$ is idempotent, so by Proposition~\ref{prop:free-implies-cofibrant}, there exists a computad $D$ and morphisms
  $\iota : D\to C$ and $\pi : C\to D$ such that
  \begin{align*}
    \pi\iota &= \id &\iota\pi &= sr
  \end{align*}
  It follows that $r\iota : D\to \bbX$ is an isomorphism with inverse
  $\pi s : \bbX\to D$, so $\bbX$ is free on a computad.
\end{proof}

\section{Higher categories}\label{sec:examples}

Over the past decades, many definitions of higher categorical structures have
been proposed, a number of which can be described as presheaves on some category
of shapes $\clI$ that are either equipped with extra operations (algebraic
models) or satisfying certain conditions (geometric models). It is often
possible in some occasions to replace the latter with the former~\cite{nikolaus_algebraic_2011, bourke_algebraically_2020}, which can more
easily described in our setting. In this section, we will explain how Leinster's
weak $\omega$\nobreakdash-categories~\cite{leinster_higher_2004}, and algebraic
semi-simplicial Kan complexes, a model of $\infty$-groupoids, are algebras
for certain signatures.
We will also propose a signature for fully weak multiple categories, inspired by
the signature for $\omega$\nobreakdash-categories. Given that opetopic higher
categories~\cite{baez_higherdimensional_1998} are presheaves on the direct category of
opetopes~\cite{cheng_category_2003}, and fair categories are presheaves on
the direct category `fat Delta'~\cite{kock_weak_2006}, we believe that
some variant of them should also be describable as algebras of some
signature.

\paragraph{Universal algebra}
The discrete category on any set $S$ can be made
into a direct category by equipping it with a constant dimension function.
Then $S$-sorted signatures are given by an $S$-indexed family of sets
$\Sigma = (\Sigma_s)_{s\in S}$ together with an $S$-indexed family of sets
$B_f$ for every $s\in S$ and $f\in \Sigma_s$. Algebras for such a signature
are again $S$-indexed families of sets $X = (X_s)_{s\in S}$ equipped with a
function
\[
  f^\bbX : \Set^S(B_f, X) \to X_s
\]
for every sort $s\in S$ and function symbol $f\in \Sigma_s$.

Letting $S$ be a singleton, one recovers the usual notion of (infinitary)
signature of universal algebra. Letting, for example, $\Sigma = \{+,0,-\}$
consist of three elements with arities a set with two elements, no elements and
a unique element respectively, one recovers the language of group theory. Groups
can be described as algebras satisfying certain equational axioms.

Modules over a ring can similarly be expressed as algebras for some signature
satisfying certain equations. The set $S$ of sorts in this case consists of two
elements $R$ and $V$, representing the ring elements and the vectors. The set
$\Sigma_R$ contains two function symbols $+^R,\cdot^R$ of arity $y_R\amalg y_R$,
for $y_\bullet$ the Yoneda embedding, two function symbols $0^R,1^R$ with arity
the empty family, and one function symbol $-^R$ with arity $y_R$. The set
$\Sigma_V$ contains similarly function symbols $+^V,0^V,-^V$ with the obvious
arities and a function symbol $\cdot^V$ with arity $y_R\amalg y_V$. We hope that
out of those examples, it is clear how to incorporate arbitrary many-sorted
signatures in our framework.

\paragraph{Kan complexes}
A classical result in algebraic topology shows that the homotopy theory of
spaces and simplicial sets are equivalent. More recently, it was shown that one
may define a weak model structure on a semi-simplicial set, presenting the
same theory. The cofibrations of this weak model structure are monomorphisms,
generating acyclic cofibrations are the horn inclusions, and weak equivalences
are the morphisms that become homotopy equivalences of spaces after
realisation~\cite[Theorem~5.5.6]{henry_weak_2020}.
This structure can be transferred to an actual model structure on the category
of algebraically fibrant objects that we discuss
below~\cite[Example~33]{bourke_algebraically_2020}.

To set the notation, let $\Delta_+$ the category of non-empty finite ordinals
and strictly monotone functions. This category has as objects natural numbers
$[n] = \{0,\dots, n\}$ and morphisms generated by the \emph{face maps}
\begin{align*}
  \delta_i^n &: [n-1] \to [n] \\
  \delta_i^n(j) &=
  \begin{cases}
    j, &\text{if }j<i \\
    j+1, &\text{if }i\le j \\
  \end{cases}
\end{align*}
for $n>0$ and $0 \le i \le n$ under the simplicial identity
\[\delta_i^{n+1}\delta_j^n = \delta_{j+1}^{n+1}\delta_i^n\]
for $n>0$ and $0\le i\le j\le n$. It is clear that $\Delta_+$ is a direct
category with dimension function $\dim([n]) = n$. Presheaves on it are called
\emph{semi-simplicial sets}.

We will denote by $\Delta[n]$ the semi-simplicial set represented by $[n]$, and
define its \emph{boundary} $\partial \Delta[n]$ to be its semi-simplicial subset
obtained by removing the top-dimensional element $\id_{[n]}$. For every $n>0$
and $0\le k\le n$, we define the semi-simplicial horn $\Lambda^k[n]$ to be the
subset of $\partial \Delta[n]$ obtained by removing the face $\delta^n_k$. An
\emph{algebraic semi-simplicial Kan complex} $\bbX$ is a semi-simplicial set $X$
equipped with a choice of lifts for every diagram of the form
\[\begin{tikzcd}
	{\Lambda^k[n]} & X \\
	{\Delta[n]}
	\arrow["\sigma", from=1-1, to=1-2]
	\arrow[dashed, from=2-1, to=1-2]
	\arrow[hook, from=1-1, to=2-1].
\end{tikzcd}\]
The choices of lifts clearly amount to two operations
\begin{align*}
  \face^\bbX_{k,n} &: [\op{\Delta}_+,\Set](\Lambda^k[n],X) \to X_{n-1} \\
  \filler^\bbX_{k,n} &: [\op{\Delta}_+,\Set](\Lambda^k[n],X) \to X_n
\end{align*}
satisfying certain boundary conditions, that provide the image of the missing
face and of the interior of $\Delta[n]$ respectively. Morphisms of algebraic
semi-simplicial Kan complexes are morphisms of presheaves preserving the
chosen lifts.

From the description above, it is easy to extract a $\Delta_+$-sorted signature
$\Sigma_{\Kan} = (\Sigma_{\Kan,n})_{n\in \bbN}$ whose algebras are algebraic
semi-simplicial Kan complexes. The signature $\Sigma_{\Kan,n}$ has two families
of function symbols of sort $[n]$. The first one consists of the symbols
$\face_{k,n+1}$ for $0\le k \le n+1$ with arity $\Lambda^k[n+1]$ and with
boundary terms
\[
  t_{\face_{k,n+1}, \delta} = \var(\delta_k^{n+1}\circ \delta)
\]
for every non-identity morphism $\delta : [m]\to [n]$. The second one consists
of the symbols $\filler_{k,n}$ for $n>0$ and $0\le k\le n$ with arity
$\Lambda^k[n]$ and with boundary terms
\[
  t_{\filler_{k,n},\delta} = \begin{cases}
    \face_{k,n}[\id], &\text{if }\delta = \delta_k^n, \\
    \var(\delta), &\text{otherwise.}
  \end{cases}
\]
It is easy to see that those families of boundary terms satisfy the cocycle
conditions for giving a type, and that they correspond to the boundary
conditions that the operations $\face_{k,n}$ and $\filler_{k,n}$ must satisfy.

As explained above, the category of algebraic semi-simplicial Kan complexes
admits a model structure equivalent to spaces. Weak equivalences of algebraic
Kan complexes are morphisms that become homotopy equivalences after
geometrically realising the underlying semi-simplicial set, while the two
weak factorisation systems are cofibrantly generated by the set of boundary
inclusions $\partial \Delta[n] \subseteq\Delta[n]$, and by the set of horn
inclusions $\Lambda^k[n]\subseteq \Delta[n]$ respectively, seen as morphisms of
free algebraic semi-simplicial sets. In particular, the (cofibrations, trivial
fibrations) weak factorisation system coincides with the one discussed in
Section~\ref{sec:computads-are-cofibrant}, so computads are the
cofibrant complexes.

\paragraph{Globular categories}
The motivating example for this work is the theory of globular weak
$\omega$\nobreakdash-categories~\cite{batanin_monoidal_1998,leinster_higher_2004}.
This example was studied extensively in our previous work~\cite{dean_computads_2022}, where a
$\bbG$-sorted signature was provided whose algebras coincide with Leinster's
$\omega$\nobreakdash-categories. The underlying category $\bbG$ of sorts is the category of
globes with objects natural numbers and morphisms generated by
\[s_n, t_n : (n) \to (n+1)\]
under the globularity conditions
\begin{align*}
  s_{n+1}s_n &= t_{n+1}s_n &
  s_{n+1}t_n &= t_{n+1}t_n.
\end{align*}
This is clearly a direct category with dimension function the identity.
Presheaves on it are called \emph{globular sets} and they can be visualised,
due to the globularity conditions, as collections of directed disks.

\begin{figure}[b]
  \[
    \begin{tikzcd}
      \bullet & \bullet & \bullet & \bullet & \bullet & \bullet & \bullet
      \arrow[from=1-2, to=1-3]
      \arrow[""{name=0, anchor=center, inner sep=0},
        bend left=45, from=1-4, to=1-5]
      \arrow[""{name=1, anchor=center, inner sep=0},
        bend right=45, from=1-4, to=1-5]
      \arrow[""{name=2, anchor=center, inner sep=0},
        bend left=45, from=1-6, to=1-7]
      \arrow[""{name=3, anchor=center, inner sep=0},
        bend right=45, from=1-6, to=1-7]
      \arrow[shorten <=3pt, shorten >=3pt, Rightarrow, from=0, to=1]
      \arrow[""{name=4, anchor=center, inner sep=0}, bend right=45,
      shorten <=3pt, shorten >=3pt, Rightarrow, from=2, to=3]
      \arrow[""{name=5, anchor=center, inner sep=0}, bend left=45,
      shorten <=4pt, shorten >=4pt, Rightarrow, from=2, to=3]
      \arrow[Rightarrow, shorten <=2pt, shorten >=2pt, from=4, to=5]
      \arrow[shorten <=2pt, shorten >=2pt, from=4, to=5, no head]
    \end{tikzcd}
  \]
  \caption{The first four representable globular sets.}
\end{figure}
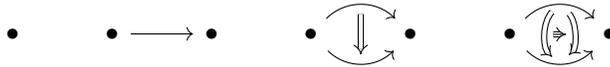

The arities for the signature $\Sigma_{\cat}$ are given by a family of globular
sets indexed by rooted planar trees, which we call \emph{Batanin trees}~\cite{batanin_monoidal_1998}. 
For alternative descriptions of this families,
see~\cite{berger_cellular_2002,leinster_higher_2004,dean_computads_2022}. This
family consists of a globular set $\Bat$ with $s_n^* = t_n^* = \partial_n$.
It consists moreover of a globular set
$\Pos B$ for every Batanin tree $B$, and a pair of inclusions
\begin{align*}
  s_n^B &: \Pos(\partial_n B) \to \Pos B &
  t_n^B &: \Pos(\partial_n B) \to \Pos B
\end{align*}
for every tree $B\in \Bat_{n+1}$ satisfying the globularity conditions.

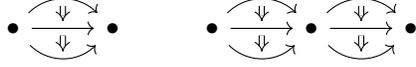
\begin{figure}[ht]
  \[
    \begin{tikzcd}
      \bullet & \bullet & \bullet & \bullet & \bullet
      \arrow[""{name=0, anchor=center, inner sep=0},
      bend left = 45, from=1-1, to=1-2]
      \arrow[""{name=1, anchor=center, inner sep=0},
      bend right = 45, from=1-1, to=1-2]
      \arrow[""{name=2, anchor=center, inner sep=0}, from=1-1, to=1-2]
      \arrow[""{name=3, anchor=center, inner sep=0},
      bend left = 45, from=1-3, to=1-4]
      \arrow[""{name=4, anchor=center, inner sep=0},
      bend left = 45, from=1-4, to=1-5]
      \arrow[""{name=5, anchor=center, inner sep=0},
      bend right = 45, from=1-3, to=1-4]
      \arrow[""{name=6, anchor=center, inner sep=0},
      bend right = 45, from=1-4, to=1-5]
      \arrow[""{name=7, anchor=center, inner sep=0}, from=1-4, to=1-5]
      \arrow[""{name=8, anchor=center, inner sep=0}, from=1-3, to=1-4]
      \arrow[shorten <=2pt, shorten >=2pt, Rightarrow, from=0, to=2]
      \arrow[shorten <=2pt, shorten >=2pt, Rightarrow, from=2, to=1]
      \arrow[shorten <=2pt, shorten >=2pt, Rightarrow, from=3, to=8]
      \arrow[shorten <=2pt, shorten >=2pt, Rightarrow, from=8, to=5]
      \arrow[shorten <=2pt, shorten >=2pt, Rightarrow, from=4, to=7]
      \arrow[shorten <=2pt, shorten >=2pt, Rightarrow, from=7, to=6]
    \end{tikzcd}
  \]
  \caption{Batanin trees corresponding to vertical composition of $2$-cells and
    the interchange axiom respectively.}
\end{figure}

This family of globular sets is well-studied, since it familially represents the
free strict $\omega$\nobreakdash-category monad. While in a strict $\omega$\nobreakdash-category,
diagrams of cells indexed by a Batanin tree admit a unique composite, in an
arbitrary $\omega$\nobreakdash-category they only admit unique composite up to a higher
coherence cell. This has be made precise in various equivalent ways using
\emph{contractible globular operads}~\cite{batanin_monoidal_1998,leinster_higher_2004}, \emph{coherators} for
$\omega$\nobreakdash-categories~\cite{maltsiniotis_grothendieck_2010}, and the type theory
\emph{CaTT}~\cite{finster_typetheoretical_2017}. Equivalences between those
approaches and the one below have already been
established~\cite{ara_infty_2010,benjamin_type_2020,dean_computads_2022}.

The signature $\Sigma_{\cat} = (\Sigma_{\cat,n})_{n\in \bbN}$ is defined
recursively. There are no function symbols of sort $0$, since there should be no
way to compose objects. The signature $\Sigma_{\cat,n+1}$ has a family of
function symbols $\coh_{B,A}$ where $B\in \Bat_{n+1}$ is a Batanin tree and
$A = (a,b)$ is a pair of terms of $\Pos B$ of sort $n$ such that
\begin{itemize}
  \item the source and target of $a$ and $b$ coincide,
  \item $a = \M_{\Sigma_{\cat,n}}(s^B_n)(a')$ for some term $a'$ of
    $\Pos(\partial_n B)$ corresponding to an epimorphism,
  \item $b = \M_{\Sigma_{\cat,n}}(t^B_n)(b')$ for some term $b'$ of
  $\Pos(\partial_n B)$ corresponding to an epimorphism.
\end{itemize}
The arity of $\coh_{B,A}$ is given by $\tr_{n+1}\Pos B$ and its top dimensional
boundary terms are given by
\begin{align*}
  t_{\coh_{B,A},s_n} &= a & t_{\coh_{B,A},t_n} &= b.
\end{align*}
The rest of the boundary terms are determined uniquely by the cocycle
conditions. The motivation for this choice of signature is that $a', b'$
represents ways to compose the boundary of $B$, which should determine a
canonical way to compose $B$. When $B\in \Bat_{n+1}$ satisfies that
$s_n^B = t_n^B = \id$, that is it has dimension at most $n$, the $a'$ and $b'$
correspond to different ways to compose $B$ and $\coh_{B,A}$ to a coherence cell
between the composites.

Similarly, we can build a signature $\Sigma_{\grpd}$ for $\omega$\nobreakdash-groupoids by
allowing as functional symbols $\coh_{B,A}$ all pairs of a Batanin tree
$B\in \Bat_{n+1}$ and pairs of terms $A = (a,b)$ with common source and target.
The motivation behind this definition being that the geometric realisation of
the globular sets $\Pos B$ is contractible, so the free $\omega$\nobreakdash-groupoid on
them should be trivially fibrant.

\paragraph{Multiple categories}
Strict multiple categories are an infinite dimensional generalisation of
$n$-fold categories. They consist of a set of objects $X_\emptyset$, a set of
arrows $X_i$ for every direction $i\in \bbN$ and more generally a set of
$n$-dimensional cubes $X_I$ for every $I\subseteq \bbN$ of cardinality $n$,
together with face maps and associative, unital composition operations
\[+_i : X_I \times_{X_{I\setminus \{i\}}} X_I\to X_I \]
for every $i\in I$ satisfying the usual interchange law. Weaker versions of
them were recently introduced~\cite{grandis_introduction_2016}, where the
composition $+_0$ is strictly associative and unital, while the rest are only
associative and unital up to a higher cell. Here, we propose an alternative
unbiased version of multiple category that is weak in all directions.

We start from the category $\bbM_+$ with objects finite subsets of the natural
numbers and a morphism
\[\delta_J^\alpha : I\setminus J\to I\]
for $J \subseteq I$ and $\alpha : J\to \{0,1\}$. Composition of morphisms is
given by disjoint union
\[\delta_J^\alpha\delta_K^\beta = \delta_{J\cup K}^{(\alpha,\beta)}\]
where $(\alpha,\beta): J\cup K\to \{0,1\}$ is the map induced by the universal
property of the disjoint union. Equivalently, morphisms in $\bbM_+$ are
generated by the face maps
\[\delta_i^\alpha : I\setminus\{i\} \to I\]
for $i\in I$ and $\alpha\in \{0,1\}$ under the commutativity condition
\[\delta_i^\alpha \delta_j^\beta = \delta_j^\beta \delta_i^\alpha.\]
Clearly, $\bbM_+$ is a direct category with dimension given by cardinality. We
will call presheaves on it, \emph{semi-multiple sets}.

Iterating the composition operations in a strict multiple category, one can
compose grid-shaped arrays of cubes in a unique manner. Since grids are
determined uniquely by their number of cubes in each direction, we let
$\Grid_I$ be the set of sequences $G : I \to \bbN$ and we define for
$J \subseteq I$,
\[
  d_J = (\delta_J^\alpha)^*: \Grid_I\to \Grid_{I\setminus J}
\]
the function forgetting the values of a sequence at $J$. This assignment easily
defined a semi-multiple set.

\begin{figure}[t]
  \[
    \begin{tikzcd}[column sep = large, row sep = large]
      \bullet & \bullet & \bullet & \bullet & \bullet \\
      \bullet & \bullet & \bullet & \bullet & \bullet
      \arrow["{(0,0)}"', from=2-1, to=2-2]
      \arrow["{(1,0)}"', from=2-2, to=2-3]
      \arrow["{(2,0)}"', from=2-3, to=2-4]
      \arrow["{(2,1)}", from=1-3, to=1-4]
      \arrow["{(1,1)}", from=1-2, to=1-3]
      \arrow["{(0,1)}", from=1-1, to=1-2]
      \arrow["{(0,0)}"{description}, from=2-1, to=1-1]
      \arrow["{(1,0)}"{description}, from=2-2, to=1-2]
      \arrow["{(2,0)}"{description}, from=2-3, to=1-3]
      \arrow["{(3,0)}"{description}, from=2-4, to=1-4]
      \arrow["{(3,0)}"', from=2-4, to=2-5]
      \arrow["{(3,1)}", from=1-4, to=1-5]
      \arrow["{(4,0)}"{description}, from=2-5, to=1-5]
      \arrow["{(0,0)}"{description}, shorten <=4pt, shorten >=4pt,
      Rightarrow, from=2-1, to=1-2]
      \arrow["{(1,0)}"{description}, shorten <=4pt, shorten >=4pt,
      Rightarrow, from=2-2, to=1-3]
      \arrow["{(2,0)}"{description}, shorten <=4pt, shorten >=4pt,
      Rightarrow, from=2-3, to=1-4]
      \arrow["{(3,0)}"{description}, shorten <=4pt, shorten >=4pt,
      Rightarrow, from=2-4, to=1-5]
    \end{tikzcd}
  \]
  \caption{The semi-multiple set of positions of the grid
    $G = (4,1,0,0,\dots)$.}
\end{figure}
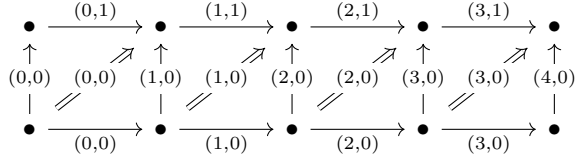

The semi-multiple set $\Pos G$ of positions of a grid $G : I\to \bbN$ consists
of the cubes of that grid. Identifying each cube with its vertex closest to the
origin, we obtain the following description. The set $\Pos_J G$ is empty unless
$J\subseteq I$, in which case it consists of points $R\in \bbN^I$ such that
$R(i)\le G(i)$ for all $i\in I$ and the inequality holds strictly when $i\in J$.
The face $d_K^\alpha R \in \Pos_{J\setminus K}G$ is given by
\[
  (d_K^\alpha R)(i) = \begin{cases}
    R(i) + \alpha(i), &\text{ when } i\in K \\
    R(i), &\text{ when } i\in I\setminus K
  \end{cases}
\]
Moreover, for every grid $G\in \Grid_I$, subset $J\subseteq I$ and morphism
$\delta^\alpha_J : I\setminus J\to I$, there exists an inclusion
\begin{align*}
  \delta_J^{\alpha,G} &: \Pos(d_J G) \to \Pos G \\
  (\delta_J^{\alpha,G} R)(i) &= \begin{cases}
    R(i), &\text{when } i\in I\setminus J \\
    \alpha(i)\cdot G(i), &\text{when }i\in J
  \end{cases}
\end{align*}
Those data forms precisely a family of semi-multiple sets.

The signature $\Sigma_{\mcat} = (\Sigma_{\mcat,n})_{n\in \bbN}$ is defined
recursively, similar to $\Sigma_{\cat}$. The signature $\Sigma_{\mcat,0}$
contains no function symbols of sort $\emptyset$. The signature
$\Sigma_{\mcat,n}$ for $n>0$ has one family of function symbols of arity $I$
for every subset $I\subset \bbN$ of size $n$. The family consists of symbols
$\coh_{G,A}$, where $G\in \Grid_I$ is a grid and
$A = (t_i^\alpha)_{i\in I, \alpha\in \{0,1\}}$ is a collection of
$\Sigma_{\mcat,n}$-terms of $\Pos G$ of sort $I\setminus\{i\}$ such that
\begin{itemize}
  \item $(\delta_j^\beta)^* t_i^\alpha = (\delta_i^\alpha)^* t_j^\beta$,
  \item $t_i^\alpha =
    \M_{\Sigma_{\mcat,n}}(\delta_i^{\alpha,G})(t_i^{\alpha\prime})$ for some
    term $t_i^{\alpha\prime}$ corresponding to an epimorphism.
\end{itemize}
The arity of $\coh_{G,A}$ is given by $\tr_n \Pos G$ and its top dimensional
boundary terms are given by
\[t_{\coh_{G,A}, \delta_i^\alpha} = t_i^\alpha.\]
The rest of its boundary terms are determined by the cocycle conditions.

\appendix

\section{A factorisation system for computads}
\label{sec:image-factorisation}

As shown in Theorem~\ref{thm:compv-topos}, the category of
computads and variable-to-variable morphisms is a presheaf topos.
As a consequence, every variable-to-variable morphism admits a factorisation
into an epimorphism followed by a monomorphism. It turns out that this is true
more generally in the category of computads and arbitrary morphisms:
epimorphisms and variable-to-variable monomorphisms form an orthogonal
factorisation system.

The inclusion $\zeta_\Sigma : \Compv_\Sigma\hookrightarrow\Comp_\Sigma$ is
faithful and preserves connected limits and colimits, so it preserves and
reflects both epimorphisms and monomorphisms. Using that the nerve functor of
Theorem~\ref{thm:compv-topos} is an equivalence, we see that a
variable-to-variable morphism $\sigma : C\to D$ is a monomorphism of computads
exactly
when the functions $\sigma_i : V^C_i\to V^D_i$ are injective for every sort~$i$.

In order to characterise epimorphisms of computads, we introduce the
\emph{support} of a morphism of computads. Intuitively the support captures the
generators of the target that are used in the definition of the morphism. We
will see that a morphism is an epimorphism exactly when its support contains
every generator of its target.

\begin{definition}\label{def:support-term}
  Let $C$ a $\Sigma$\nobreakdash-computad and $i$ a sort. The \emph{support of sort $i$} of
  a term $t\in \Term_{\Sigma,j}(C)$ of some sort $j$ of dimension
  $\beta \le \alpha$ is defined recursively by
  \begin{align*}
    \supp_i(\var v)
      &= \{ v \;:\; i = j\}\cup
      \bigcup_{\substack{\id\not =\delta : k \to j}} \supp_i(\phi_\delta^C(v)) \\
    \supp_i(f[\tau])
      &= \bigcup_{\substack{k\in \clI_{\beta} \\b\in B_{f,k}}}\supp_i(\tau_k(b))
  \end{align*}
  The \emph{support of sort $i$} of a morphism $\sigma : D\to C$ is defined
  by
  \[
    \supp_i(\sigma) =
    \bigcup_{\substack{j\in \clI_\alpha \\ v\in V_j^D}} \supp_i(\sigma_j(v)).
  \]
\end{definition}

Recursively, we can see that the support of sort $i$ of a term of sort
$j$ is empty unless there exists a morphism $\delta : j\to i$, which in
particular implies that $\dim j \le \dim i$. The following lemma shows that the
support is closed under boundary maps, which implies in particular that the
support of a term and the morphism classifying it coincide. Its proof
is a simple mutual induction, left to the reader.

\begin{lemma}\label{lem:support-boundary}
  Let $C$ a computad, $i$ a sort, and $t\in \Term_{\Sigma,j}(C)$ a
  term of $C$ of some sort $j$ of dimension $\beta\le\alpha$.
  For every face map $\delta : k\to j$,
  \[\supp_i(\delta^*t)\subseteq \supp_i(t),\]
  while for every morphism of $\Sigma$\nobreakdash-computads $\sigma : C\to D$,
  \[
    \supp_i(\Term_\Sigma(\sigma)(t)) =
    \bigcup_{\substack{k\in \clI_{\beta} \\ v\in \supp_k(t)}}
    \supp_i(\sigma_k(v)).
  \]
\end{lemma}

\begin{proposition}\label{prop:support-equalises}
  Two morphisms of computads agree on a term if and
  only if they agree on every generator of its support.
\end{proposition}
\begin{proof}
  Let $\sigma, \sigma' :C\to D$ morphisms of computads. By induction on
  the dimension of the signature, it suffices to show the result for
  terms $t$ of $C$ of sort $j$ of dimension $\alpha$. Suppose first
  that $t = \var v$ is a generator. Since $v\in \supp_j(v)$, one direction
  obviously holds. For the converse, suppose that $\sigma$ and $\sigma'$ agree
  on $t$. Then they must also agree on $\phi_\delta^C(v)$ for every
  $\delta : k\to j$. By the inductive hypothesis, they agree on the support of
  all $\phi_\delta^C(v)$, hence on the support of $t$ as well.

  Suppose now that $t = f[\tau]$. Then $\sigma$ and $\sigma'$ agree on $t$
  if and only if $\sigma\tau = \sigma\tau'$. This amounts to them agreeing on
  $\tau_k(b)$ for every $k\in \clI_{\dim j}$ and $b\in B_{f,k}$. By the
  inductive hypothesis, this is equivalent to agreeing on the support of each
  $\tau_k(b)$. Equivalently, they agree on the union of those supports, which
  is the support of $t$.
\end{proof}

\begin{corollary}\label{cor:full-implies-epi}
  If the support of a morphism contains every generator, then it is an
  epimorphism.
\end{corollary}

The converse of this corollary also holds. We will deduce it by showing that
every morphism of computads factors uniquely as one whose support contains all
generators followed by one a variable-to-variable monomorphism. the following
lemma concerning the support of a variable-to-variable morphism can easily be
shown again by induction on depth.

\begin{lemma}\label{lem:supp-var-to-var}
  The support of a variable-to-variable morphism $\sigma : C\to D$ consists
  precisely of the generators in the images of the functions $\sigma_i$.
\end{lemma}

\begin{lemma}\label{lem:supp-determines-lifts}
  Let $\rho : D\to E$ a variable-to-variable monomorphism and let
  $\sigma : C\to D$ arbitrary morphism. Then $\sigma = \rho\sigma'$ for some
  $\sigma' : C\to D$ if and only if the support of $\sigma$ is contained in
  that of $\rho$. Moreover, the factorisation is unique.
\end{lemma}

\begin{proof}
  Uniqueness of the factorisation follows by $\rho$ being a monomorphism, while
  necessity of this condition for the existence of a factorisation follows by
  Lemma~\ref{lem:support-boundary}. By induction on the dimension of the
  signature, we may assume that this condition is sufficient for morphisms of
  $\Sigma_\beta$-computads for all $\beta < \alpha$.

  We proceed by induction on the depth of $\sigma$. Suppose first that
  $C = \bbD^i_\Sigma$ and that $\sigma$ classifies a generator $t = \var v$. By
  Lemma~\ref{lem:supp-var-to-var}, there exists $v'\in V_i^D$, mapped via
  $\rho$ to $t$. The morphism corresponding to $\var v'$ gives the claimed
  factorisation. If $\sigma$ classifies a term $t = f[\tau]$, then
  \[\supp_i(\sigma) = \supp_i(t) = \supp_i(\tau),\]
  so by the inductive hypothesis, there exists a factorisation
  $\tau = \rho\tau'$. The morphism corresponding to $f[\tau']$ gives a
  factorisation of $\sigma$ via $\rho$.

  Let now $C$ and $\sigma$ be arbitrary. By the inductive hypothesis, there
  exist unique $\sigma'_\beta : C_\beta\to D_\beta$ such that
  $\sigma_\beta = \rho_\beta\sigma_\beta'$ for all $\beta < \alpha$ and unique
  terms $\sigma'_i(v)\in \Term_\Sigma(D)$ for every generator $v\in V_i^C$ of
  dimension $\alpha$ such that
  \[\Term_\Sigma(\rho)(\sigma'_i(v)) = \sigma_i(v).\]
  By uniqueness, this data assembles into a morphism $\sigma' : C\to D$ such
  that $\sigma= \rho\sigma'$.
\end{proof}

\begin{proposition}\label{prop:full-mono-factorisation}
  Every morphism can be factored as a morphism whose support contains all
  generators, followed by a variable-to-variable monomorphism.
\end{proposition}

\begin{proof}
  Let $\sigma : C\to D$ arbitrary morphism. We will construct recursively a
  variable-to-variable monomorphism $\iota_\sigma : \supp\sigma\to D$ where the
  generators of $\supp\sigma$ are given by $V_i^{\supp\sigma} = \supp_i\sigma$,
  and $\sigma_i$ is the obvious subset inclusion. The gluing maps of
  $\supp\sigma$ are defined recursively: given a face map $\delta : j\to i$,
  and some $v\in \supp_i\sigma$, by Lemma~\ref{lem:supp-determines-lifts},
  there exists unique term $\phi_\delta^{\supp\sigma}(v)$ of
  $\supp\sigma$ satisfying that
  \[
    \Term_{\Sigma_{\beta}}(\iota_{\sigma,\beta})\phi_\delta^{\supp\sigma}(v)
    = \phi_\delta^{C}(v),
  \]
  where $\beta = \dim j$. Uniqueness implies the cocycle condition
  $\delta^*\phi_{\delta'}^{\supp\sigma} = \phi_{\delta'\delta}^{\supp\sigma}$,
  so this data defines a computad $\supp\sigma$ and a
  variable-to-variable monomorphism out of it. Moreover, by construction,
  \[\supp_i(\iota_\sigma) = \supp_i(\sigma),\]
  so by the same lemma, there exists unique
  \[\pi_\sigma : C\to \supp\sigma\]
  such that $\sigma = \pi_\sigma\iota_\sigma$. A simple recursive argument
  shows that the support of $\pi_\sigma$ contains all generators of
  $\supp\sigma$.
\end{proof}

\begin{corollary}\label{cor:epi-implpies-full}
  The support of epimorphisms contains all generators.
\end{corollary}
\begin{proof}
  If $\sigma$ is an epimorphism, then the variable-to-variable monomorphism
  $\iota_\sigma$ is also epic. Since $\Compv_\Sigma$ is a topos,
  $\iota_\sigma$ must be invertible, hence bijective on generators. Therefore,
  its support contains all generators. Such morphisms are closed under
  composition by Lemma~\ref{lem:support-boundary}, so the support of
  $\sigma = \pi_\sigma\iota_\sigma$ must contain all generators as well.
\end{proof}

\begin{corollary}\label{cor:epi-mono-ofs}
  Epimorphisms and variable-to-variable monomorphisms form an orthogonal
  factorisation system.
\end{corollary}
\begin{proof}
  Both classes contain isomorphisms and they are closed under composition.
  Moreover, every morphism factors as an epimorphism followed by a
  variable-to-variable
  monomorphism, so it remains to show tha the factorisation is unique
  up to unique isomorphism~\cite[Proposition~14.7]{adamek_abstract_1990}.
  Uniqueness follows from the left class being epimorphism, so it remains to
  show existence. For that, let $\sigma : C\to D$ an epimorphism,
  $\rho : D\to E$ a variable-to-variable monomorphism and consider the
  commutative square
  \[\begin{tikzcd}
    C
    \rar{\pi_{\rho\sigma}}
    \dar[swap]{\sigma} &
    \supp\sigma
    \dar{\iota_\sigma} \\
    D
    \urar[dashed]{\chi}
    \rar{\rho} &
    E
  \end{tikzcd}\]
  Since $\sigma$ is epic, the supports of $\rho$, $\sigma\rho$ and
  $\iota_{\sigma\rho}$ coincide. Two applications of Lemma~\ref{lem:supp-determines-lifts} give a diagonal lift $\chi$ as in the diagram
  and a lift $\chi^{-1}$ in the opposite direction. By definition,
   \begin{align*}
    \iota_{\rho\sigma}\chi\chi^{-1} &= \iota_{\rho\sigma} &
    \rho\chi^{-1}\chi &= \rho.
  \end{align*}
  Since $\rho$ and $\iota_{\rho\sigma}$ are monic, the morphism $\chi$ is an
  isomorphism with inverse $\chi^{-1}$.
\end{proof}

\begin{corollary}\label{cor:comp-cauchy}
  The category of computads is Cauchy-complete.
\end{corollary}
\begin{proof}
  Let $\sigma: C\to C$ an idempotent morphism of computads. By Corollary~\ref{cor:epi-mono-ofs}, there exists a factorisation $\sigma = \iota\pi$
  with $\pi$ epic and $\iota$ monic. Since $\sigma\sigma = \sigma$, it follows
  that $\pi\iota = \id$. Therefore, every idempotent morphism of computads
  splits.
\end{proof}

\section{Some properties of the term monad}\label{sec:properties-term-monad}

We conclude with some technical properties of the term monad, namely
that it is cartesian and accessible. From the latter, we will deduce that the
category of $\Sigma$\nobreakdash-algebras is locally presentable, hence complete,
cocomplete, and every set of morphisms cofibrantly generates a weak
factorisation system. As usually, $\Sigma$ will denote an $\clI$-sorted
signature of some dimension $\alpha$.

\begin{proposition}\label{prop:term-monad-cartesian}
  The term monad is cartesian.
\end{proposition}

\begin{proof}
  Let $\Cptdv_\Sigma : [\op{\clI_\alpha},\Set]\to \Compv_\Sigma$ the
  restriction of the inclusion of presheaves into computads to the subcategory
  of variable-to-variable morphisms. The composite
  $V_i^\bullet\circ \Cptdv_\Sigma$ is cocontinuous for every sort $i$, since it
  sends a presheaf $X$ to the set $X_i$. By the decomposition of $V_i^\bullet$
  in Section~\ref{sec:computads-are-presheaves}, the functor $\clT_p\circ
  \Cptdv_\Sigma$ preserve connected limits for every plex $p$. By
  Theorem~\ref{thm:compv-topos}, $\Cptdv_\Sigma$ must preserve them as well.
  Corollary~\ref{cor:termv-fam-representable} finally implies that the composite
  $\M_\Sigma = \Termv_\Sigma\Cptdv_\Sigma$ must preserve them as well.

  It remains to show that $\eta_\Sigma$ and $\mu_\Sigma$ are cartesian natural
  transformations. For the former, let $\sigma : X\to Y$ a morphism of
  presheaves and consider the naturality square
  \[
    \begin{tikzcd}
      X
      \rar{\sigma}
      \dar{\eta_{\Sigma,X}} &
      Y \dar{\eta_{\Sigma,Y}} \\
      \M_\Sigma(X)
      \rar{\M_\Sigma(\sigma)} &
      \M_\Sigma(Y)
    \end{tikzcd}
  \]
  Pullbacks in categories of presheaves are computed object-wise, so we need to
  show for every sort $i$ that pairs $t\in \M_{\Sigma,i}(X)$ and $y\in Y_i$ that
  satisfy the compatibility condition
  \[ \M_\Sigma(\sigma)(t) = \eta_{\Sigma,Y}(y) = \var y \]
  can be lifted uniquely to an element of $X_i$. From the compatibility
  condition, we deduce that $t$ must be a generator, so $t = \eta_{\Sigma,X}(x)$
  for unique $x\in X_i$. Substituting $t$ into the compatibility condition, we
  get then that $y = \sigma(x)$, so the square is a pullback.

  The multiplication of the monad can be written as the following whiskered
  composite
  \[
    \mu_\Sigma =
    \Term_\Sigma\varepsilon_\Sigma(\zeta_\Sigma\Cptd_\Sigma^{\var}).
  \]
  The functor of $\Sigma$\nobreakdash-terms is representable by $\bbD^i_\Sigma$, so it
  preserves pullbacks. Therefore, it suffices to show that
  $\varepsilon_\Sigma\zeta_\Sigma$ is cartesian. Let therefore $\rho : C\to D$
  a variable-to-variable morphism of computads and suppose that a
  solid commutative diagram of the following form is given
  \[
    \begin{tikzcd}[column sep = huge]
      E
      \drar[dashed]{\tau}
      \arrow["\tau^1", from =1-1, to=2-3, bend left = 15]
      \arrow["\tau^2"', from =1-1, to=3-2, bend right]
      \\
      &
      \Cptd_\Sigma\Term_\Sigma C
      \rar{\Cptd_\Sigma\Term_\Sigma\rho}
      \dar{\varepsilon_{\Sigma,C}} &
      \Cptd_\Sigma\Term_\Sigma D
      \dar{\varepsilon_{\Sigma,D}} \\
      &
      C
      \rar{\rho} &
      D
    \end{tikzcd}
  \]
  By induction on the dimension of the signature, we may assume that
  $\varepsilon_{\Sigma_\beta}\zeta_{\Sigma_\beta}$ is cartesian for all
  $\beta < \alpha$. We proceed by recursion on the depth of $\tau^1$ to show
  that there exists unique $\tau$ making the entire diagram above commute.

  Suppose first that $E = \bbD^i_\Sigma$ for some sort $i$ of sort $\alpha$, and
  suppose further that $\tau^1$ classifies a generator $t_1 = \var t_1'$ where
  $t_1'\in \Term_{\Sigma,i}(D)$. Then $\tau^2$ classifies some term
  $t_2\in \Term_{\Sigma,i}(C)$ such that
  \[\Term_\Sigma(\rho)(t_2) = t_1'.\]
  The morphism $\tau$ must classify some generator for the lower triangle to
  commute, and that generator must by $\var(t_2)$ for the left triangle to
  commute. Conversely, the morphism classifying $\var t_2$ makes the diagram
  commute, so there exist unique $\tau$ making the diagram commute.

  Suppose then that $\tau^1$ classifies a composite term
  $t_1 = f[\hat{\tau}^1]$ and let $t_2$ the term classified by $\tau^2$.
  Then
  \[\Term_\Sigma(\rho)(t_2) = f[\varepsilon_{\Sigma,D}\hat{\tau}^1]\]
  and $\rho$ is variable-to-variable, so $t_2 = f[\hat{\tau}^2]$ for
  some $\hat{\tau}^2$ satisfying that
  \[\rho\hat{\tau}^2 = \varepsilon_{\Sigma,D}\hat{\tau}^1.\]
  By the recursive hypothesis, there exists unique
  $\hat{\tau} : B_f\to \Cptd_\Sigma\Term_\Sigma C$ such that
  \begin{align*}
    \hat{\tau}^1 &= (\Cptd_\Sigma\Term_\Sigma\rho)\hat{\tau} &
    \hat{\tau}^2 &= \varepsilon_{\Sigma,C}\hat{\tau}.
  \end{align*}
  The morphism $\tau$ corresponding to the term $f[\hat{\tau}]$ makes the
  diagram above commute, and it is easily seen to be unique.

  Finally, let $E$ be arbitrary. Then by the recursive hypothesis, for every
  $\beta <\alpha$, there exists unique
  $\tau_\beta : E_\beta\to \tr_\beta^\Sigma\Cptd_\Sigma\Term_\Sigma C$
  making the obvious truncated versions of the diagram above and there exist
  for every $v\in V_i^E$ for $i$ of dimension $\alpha$, unique term $\tau_i(v)$
  of $\Cptd_\Sigma\Term_\Sigma C$ such that
  \begin{align*}
    \Term_\Sigma(\varepsilon_{\Sigma,C})(\tau_i(v)) &= \tau_i^2(v) &
    (M_\Sigma\Term_\Sigma\rho)(\tau_i(v)) &= \tau^1_i(v).
  \end{align*}
  By uniqueness, we can easily deduce that those morphisms and terms constitute
  a morphism $\tau : C\to D$ making the diagram above commute. It is not hard to
  see that said $\tau$ is unique.
\end{proof}

\begin{proposition}\label{prop:term-monad-accessible}
  The monad $\M_\Sigma$ preserves $\lambda$-filtered colimits for some regular
  cardinal $\lambda$.
\end{proposition}

\begin{proof}
  Let $\lambda$ a regular cardinal strictly greater than the cardinality of the
  disjoint union $\coprod_{j} B_{f,j}$ for every sort $i\in \clI_\alpha$, and
  every function symbol $f\in \Sigma_{i,F}$. Such $\lambda$ exists, since both
  the category of sorts and the collection of function symbols are small. To
  show that $\M_{\Sigma}$ preserves
  $\lambda$-filtered colimits, it suffices to show that $\M_{\Sigma,i}$
  preserves them for every sort $i$. By induction on the dimension $\alpha$ of
  the signature, we may assume that $\M_{\Sigma_\beta}$ preserves
  $\lambda$-filtered colimits for all $\beta < \alpha$, so it remains to show
  that $\M_{\Sigma,i}$ preserves them for $i$ of dimension $\alpha$.

  In order to do that, define for every ordinal $\gamma \le \lambda$,
  \[\M_{\Sigma,i}^{\gamma} : [\op{\clI_\alpha},\Set] \to \Set\]
  to be the functor sending a presheaf $X$ to the set of terms of
  $\Cptd_\Sigma X$ of sort $i$ and recursive depth at most $\gamma$. Define also
  for every function symbol $f\in \Sigma_i$, a functor
  \[\M_{\Sigma,i,f}^{\gamma} : [\op{\clI_\alpha},\Set]\to \Set\]
  sending a presheaf $X$ to the set of morphisms
  $\Cptd_\Sigma(B_f)\to \Cptd_\Sigma(X)$ of depth at most $\beta$. The
  discussion on recursive depth shows that for $\beta = \lambda$, we recover the
  sets of all terms of sort $i$, and all morphisms respectively, so it suffices
  to show that those functors preserve $\lambda$-filtered colimits.

  We proceed recursively on the ordinal $\gamma$. The functor
  $\M_{\Sigma,i}^{0}$ is cocontinuous, since it sends $X$ to the set $X_i$.
  The functor $\M_{\Sigma,i}^{\gamma+1}$ for $\gamma < \lambda$ decomposes as
  \[
    \M_{\Sigma,i}^{\gamma+1} = \M_{\Sigma,i}^{0} \amalg
    \coprod_{f\in \Sigma_i} \M_{\Sigma,i,f}^{\gamma}
  \]
  while for $\gamma\le\lambda$ limit ordinal, we have that
  \[
    \M_{\Sigma,i}^{\gamma} =
    \colim_{\gamma'<\gamma}\M_{\Sigma,i}^{\gamma'}.
  \]
  In both cases, we see that the functor preserves $\lambda$-filtered
  colimits by the inductive hypothesis and commutativity of colimits with
  colimits. Finally, given any function symbol $f\in \Sigma_i$ and
  $\gamma\le\lambda$ arbitrary, consider the functor from the category of
  elements of the presheaf $B_f$ sending $j\in \clI_\alpha$ and $b\in B_{f,j}$
  to the set $\M_{\Sigma,j}$ when $\dim j < \alpha$, and to
  $\M_{\Sigma,j}^{\gamma}$ otherwise. The domain of the functor is
  $\lambda$-small and its limit is $\M_{\Sigma,i,f}^{\gamma}$. By the
  inductive hypothesis and commutativity of $\lambda$-small limits with
  $\lambda$-filtered colimits, we see that $\M_{\Sigma,i,f}^{\gamma}$
  preserves $\lambda$-filtered colimits. This concludes the induction.
\end{proof}

\begin{corollary}\label{cor:algebras-presentable}
  The category of algebras is locally presentable.
\end{corollary}
\begin{proof}
  Let $\lambda$ a regular cardinal such that $\M_\Sigma$ preserves
  $\lambda$-filtered colimits. The category $[\op{\clI_\alpha},\Set]$ is locally
  finitely presentable, being a presheaf
  category~\cite[Example~1.1.12]{adamek_locally_1994}, hence also
  $\lambda$-presentable~\cite[Remark~1.1.20]{adamek_locally_1994}.
  The category of algebras
  is therefore the category of algebras of a $\lambda$-accessible monad on a
  locally $\lambda$-presentable category. Hence, it is also a locally
  $\lambda$-presentable category~\cite[Remark~2.2.78]{adamek_locally_1994}.
\end{proof}

\bibliography{bibliography}
\bibliographystyle{plainurl}

\end{document}